\documentclass[11pt]{article}
\usepackage{algorithm}
 \usepackage{algorithmicx}
\usepackage{algpseudocode}

\usepackage{amsfonts}
\usepackage{bm}
\usepackage{amsmath,amssymb,amsthm}
\usepackage{enumitem}
\usepackage{natbib}
\usepackage{ulem}   % for highlighting and emphasis
\usepackage{xcolor} % use colors in text and equations
\usepackage{graphicx} %use graph format
\usepackage{epstopdf}
\usepackage{multirow}
\usepackage{subfigure}
\usepackage{caption}
\usepackage{algorithm}
\usepackage{ulem}
\usepackage{latexsym,lscape}

\usepackage[top=1in,bottom=1in,left=1in,right=1in, a4paper
]{geometry}  
\usepackage{setspace}

\newcommand{\real}{\mathbb{R}}

\allowdisplaybreaks

\makeatletter
\def\BState{\State\hskip-\ALG@thistlm}
\makeatother

\floatname{algorithm}{Procedure}

\title{Review on Ranking and Selection: A New Perspective\footnote{This paper has been accepted by {\it Frontiers of Engineering Management}.}}

\author{
        L. Jeff Hong\vspace{-8pt}\\
        {\footnotesize
        School of Management and School of Data Science, Fudan University, Shanghai 200433, China }\vspace{-8pt}\\ {\footnotesize hong\_liu@fudan.edu.cn}\\
        Weiwei Fan\thanks{Weiwei Fan is the corresponding author. Email: wfan@tongji.edu.cn.
Tel: 0086-21-65981165.} \vspace{-8pt}\\
        {\footnotesize
        Advanced Institute of Business, Tongji University, Shanghai 200092, China }\vspace{-8pt}\\ {\footnotesize wfan@tongji.edu.cn }\\
                Jun Luo\thanks{Jun Luo is the corresponding author.  Email: jluo\_ms@sjtu.edu.cn.
Tel: 0086-21-52301579. } \vspace{-8pt} \\
{\footnotesize Antai College of Economics and Management, Shanghai Jiao Tong University, Shanghai 200030, China} \vspace{-8pt} \\ {\footnotesize  jluo\_ms@sjtu.edu.cn}
}
\date{\today}

\begin{document}
\doublespacing
\maketitle

\begin{abstract}
In this paper, we briefly review the development of ranking and selection (R\&S) in the past 70 years, especially the theoretical achievements and practical applications in the past 20 years. Different from the frequentist and Bayesian classifications adopted by \cite{KimNelson2006b} and \cite{Chick2006} in their review articles, we categorize existing R\&S procedures into fixed-precision and fixed-budget procedures as in \cite{HunterNelson2017}. We show that these two categories of procedures essentially differ in the underlying methodological formulations, i.e., they are built on hypothesis-testing and dynamic-programming, respectively. In light of this variation, we review in detail some well-known procedures in the literature, and show how they fit into these two formulations. In addition, we discuss the use of R\&S procedures in solving various practical problems, and propose what we think are the important research questions in the field.

\end{abstract}

\emph{Keywords:} ranking and selection, hypothesis testing, dynamic programming, simulation

\section{Introduction}

Decision-making processes often involve comparisons among a set of alternatives regarding certain performance measures. In this study, we consider such comparison problems with the goal of selecting the best alternative, where the best is defined to have the largest (or smallest) mean performance. This is not trivial in the stochastic environment where the mean performances of these alternatives are unknown and have to be inferred via statistical sampling from stochastic systems. Therefore, a selection procedure is required to determine how many samples need to be collected from each alternative and then which alternative should be selected as the best based on the sample information. Such selection problems are often called ranking and selection (R\&S) in the literature.

R\&S problems date back to the 1950s in agricultural and clinical applications \citep{bechhofer1954single, gupta1956on}. At that time, testing the homogeneity of multiple alternatives was common (e.g., grain yields and drug treatments). For instance, an individual might desire to test whether multiple grains produced the same mean yield, or whether multiple drug treatments led to the same mean efficacy. Once the homogeneity of their means was rejected statistically, a natural issue readily arose, that is, which one was the best. {This issue was first proposed by \cite{10.2307/2236807} and triggered the early developments of R\&S.}

In the 1950s, samples needed to be collected through physical experiments, e.g., agricultural experiments and clinical trials, which might cost a long time to conduct. Thus, the experiments were often conducted in batches. Accordingly,  a considerable number of the R\&S procedures designed then were stage-wise, where the best one was selected at the end of the last stage. Starting in the 1990s, this paradigm began to change owing to the increasing computing power. An increasing number of experiments were conducted in computer simulation environments because it cost little time to generate samples. Through these simulations, samples were often collected sequentially, especially when the program is executed in a single-processor environment. This sequential nature boosted the development of sequential R\&S procedures. Unlike stage-wise procedures, sequential procedures typically provide a decision rule at each time of the sample collection process, and therefore are more efficient in most situations by taking advantage of the interim sample information. Sequential R\&S procedures are still prevalent today.

In recent years, another forming paradigm that considers large-scale R\&S problems has emerged. For early applications such as agricultural problems, the number of alternatives is relatively small. Designed for these applications, classic procedures were typically applied to problems with  fewer than 500 alternatives. However, in the modern world, we often face problems that may have thousands to tens of thousands of alternatives. For instance, in scheduling problems, one may need to determine multiple components simultaneously, such as the jobs to be scheduled, the values assigned to the jobs and the time when the scheduling happens. Assuming that 50 choices are available for each component, their combination fairly leads to a total of 125000 alternatives, which is a huge number for classic R\&S procedures. Recently, research on how to handle large-scale R\&S problems {has drawn significant attention}. As pioneer works, {\cite{LuoHong2011}, \cite{LHNW2015}} and {\cite{NHH2013, NCHH2017}} addressed large-scale problems by adapting the classic procedures into parallel computing environments.

{Interested readers may refer to \cite{FuHenderson2017} for an interesting introduction on the history of R\&S.} Basically, R\&S procedures provide a general tool for solving selection problems. Therefore they are widely applicable to practical problems. Besides, many of the R\&S procedures are also easy to implement, and some of them have been embedded in commercial simulation software packages, such as Arena and Simio.

To organize the R\&S procedures, existing review articles often categorize them into frequentist and Bayesian procedures, according to the probability models used to describe the collected samples (Chick, 2006; Kim and Nelson, 2006b; Branke et al., 2007). In this work, we take a different perspective and categorize them into fixed-precision and fixed-budget procedures, as in {\cite{NIPS2012_4640}} and \cite{HunterNelson2017}. Particularly, fixed-precision procedures intend to provide a desired statistical guarantee of the selected alternative being the best (or at least close to the best), while fixed-budget procedures intend to allocate a given sampling budget in various optimal or approximately optimal ways. To explain these two categories of procedures, we show they essentially follow two different formulations, i.e., the hypothesis-testing and dynamic-programming formulations, respectively. {A number of studies in the literature have adopted the same perspective and designed new procedures under the two formulations, e.g., \cite{BATUR2012661} and \cite{8267207}. Different from these works, the goal of this review is to construct a unified framework for each formulation and explain how the existing procedures fit in the framework.}

This paper only focuses on selecting the best mean.However, some related problems may also be categorized into R\&S problems. They essentially have different combinations of goals to achieve and performance measures used for comparisons. For instance, the goals can be ranking all the alternatives, or selecting the top $m$ alternatives, or selecting a subset of alternatives that contains the best. Meanwhile, the performance measures used can be the quantile or proportion. These problems are not covered in this study and interested readers may refer to \cite{santner1995design}, \cite{goldsman1998comparing} and \cite{KimNelson2006b} for comprehensive reviews.

{One problem closely related to R\&S is the multi-armed bandit (MAB) problem in the machine learning literature. Both problems stemmed from \cite{bechhofer1954single} and \cite{paulson1964sequential}, and they have grown into two branches of research with different goals in designing procedures. R\&S procedures typically attempt to optimize the quality of final selection. In contrast, MAB procedures attempt to balance the tradeoff between exploration (gathering new information on different alternatives) and exploitation (choosing the best alternative) in the sequential sampling process. Therefore, the MAB problem often aims to minimize cumulative regret during the sampling process. Nonetheless, a series of works have considered the pure-exploration version of the MAB problem, which is known as the best-arm identification (BAI) problem \citep{bubeck2012regret}. Although BAI and R\&S problems have the same goal, they typically make different assumptions on the samples from alternatives. Particularly, the BAI problem assumes the samples to be bounded or sub-Gaussian distributed, whereas the R\&S problem typically assumes they are Gaussian distributed with unknown variances. In this study we will not review MAB procedures. Interested readers may refer to \cite{Even-Dar2002}, \cite{bubeck2012regret}, \cite{NIPS2012_4640}, and \cite{pmlr-v30-Kaufmann13} for more information on MAB, and \cite{ma2017an} and \cite{glynn2018selecting} for their connections to R\&S procedures.}

The rest of paper is organized as follows. In Section 2, we provide a comprehensive description on how R\&S problems under fixed-precision and fixed-budget are formulated as hypothesis-testing or dynamic-programming problems, respectively. In Sections 3 and 4, we present several well-known fixed-precision and fixed-budget R\&S procedures, and explain how they can be derived under two different formulations, respectively. In Section 5, we present the procedures designed for solving large-scale R\&S problems. In Section 6, we introduce several emerging R\&S problems, followed by the discussion of some interesting future research directions in Section 7.

\section{Two formulations for R\&S}
Suppose there are $k\geq 2$ alternatives with mean performance $\pmb{\mu} = (\mu_1, \mu_2, \dots, \mu_k)$, and the best alternative is defined to have the largest mean. For simplicity, we assume that the best alternative is unique. The goal of R\&S is to select the index of the best alternative, which is unknown a priori. {If multiple alternatives have tied best means, choosing any of these alternatives as the best can be viewed as a correct selection.}

Evidently, the selection decision should be made based on the information collected from samples. Ideally, we hope to select the best alternative with 100\% probability. However, this is impossible unless infinite samples can be collected. Therefore, a tradeoff exists between the sampling budget and the precision of the selection decision. To alleviate this tradeoff, R\&S problems are often imposed with two constraints: Fixed precision and fixed budget \citep{HunterNelson2017}. In particular, the fixed-precision R\&S problems intend to achieve a fixed precision of selection when using as few sampling budget as possible, while the fixed-budget R\&S problems intend to optimize the precision of selection given a fixed sampling budget.

In this section, we show that these R\&S problems under the two constraints can be formulated as hypothesis-testing (HT) and dynamic-programming (DP) problems, respectively. We also illustrate some key issues in designing corresponding R\&S procedures.

%Let us first set up some notations. Let $X_{ij}$ denote the $j$th observation from alternative $i$, for $i=1,2,\dots,k$ and $j=1,2,\dots$. Unless specifically stated, we assume these observations are independent across alternatives and  $\{X_{ij}: j=1,2,\dots\}$ are i.i.d.\ normal with mean $\mu_i$ and variance $\sigma_i^2$. Let $X_i(n)$ and $S_i^2(n)$ denote the sample mean and sample variance calculated from the first $n$ samples from alternative $i$.

\subsection{Fixed-precision R\&S}\label{sec:HT}

To describe the precision of selection (i.e., the first constraint), one common way is to use the probability that the selected alternative is the true best, which is called the probability of correct selection (PCS). Then, under a fixed precision $1-\alpha\, ({0<\alpha<1-1/k})$, the goal of R\&S is to deliver a PCS guarantee as
\begin{align}\label{PCS}
\mbox{PCS}(\pmb \mu) = \mathrm{P}\big\{\mbox{Select the best alternative}\,|\,\pmb\mu\big\}\geq 1-\alpha, \ \forall\ \pmb\mu \in \Theta,
\end{align}
where $\Theta = \big\{\pmb\mu: \mu_{[k]}>\mu_{[k-1]}\big\}$ and $\mu_{[k]}>\mu_{[k-1]}\geq \dots\geq \mu_{[1]}$ denote the ordered means.
%As mentioned before, the R\&S problems follow from a hypothesis test of homogeneity which want to know whether these alternatives have the same means, i.e., $H_0: \mu_1=\mu_2=\dots=\mu_k$. Once $H_0$ is rejected, the R\&S problems come out naturally to ask which alternative has the largest mean.
\subsubsection{Fixed-precision R\&S formulated as hypothesis-testing}\label{subsec:2general}
Practically, one may select any alternative as the best, and then what is needed  is to tell whether this alternative is truly the best. This suffices to detect, for any alternative $j$, whether it has a larger mean than all the others, i.e., $\mu_{j}>\mu_i$ for any $i\neq j$. In light of this, the R\&S problems essentially involve $k$ simultaneous HTs and therefore {are} formulated as a multiple HT problem,
\begin{align}\label{test1}
(HT_j):  \quad\quad H_0^j: \mu_{j} \leq \max_{i\neq j}\mu_i, \quad \mbox{ versus } \quad H_1^j:  \mu_{j} > \max_{i\neq j}\mu_i, \quad \forall\, j=1,2,\dots,k.
\end{align}
Each single $HT_j$ above regards the comparison between alternative $j$ and all the others.

When $H_0^j$ is rejected, alternative $j$ should be selected as the best. Therefore, to select the best alternative correctly, we only need to avoid committing the Type II error for each $HT_j$. To make it clear, notice that the PCS guarantee in \eqref{PCS} can be rewritten as,
\begin{align*}
\mbox{PCS}(\pmb\mu) =\mathrm{P}\big\{\mbox{Reject }H_0^j\,|\, \pmb \mu \in H_1^j\big\}=1-\mathrm{P}\big\{\mbox{Type II Error in $HT_j$}\big\}\geq 1-\alpha, \mbox{ for } \pmb \mu\in H_1^j, \forall j.
\end{align*}
(For simplicity of the notation, we write $\pmb \mu\in H_d^j \,(d=0, 1)$, if $\pmb\mu$ satisfies the corresponding hypothesis.) This implies that we only need to control the Type II error for all $HT_j$ in \eqref{test1} as
\begin{align}\label{ErrorControl}
\mathrm{P}\big\{\mbox{Type II Error in $HT_j$}\big\}\leq \alpha,\quad \forall\, \pmb\mu\in H_1^j,\  j=1,2,\dots,k.
\end{align}

The Type I error for each $HT_j$ has been automatically controlled at the same time. Taking the special case when there are only two alternatives for example, when \eqref{test1} has two HTs, then the Type I error in one HT essentially corresponds to the Type II error in the other. For the general case, all $H_1^j (j=1,2,\dots,k)$ compose a disjoint partition of the whole mean space $\Theta$. This partition indicates that any mean vector $\pmb \mu$ satisfying $H_0^j$ must satisfy one of $H_1^l (l\neq j)$. Then, we are able to show
\begin{align*}
\mathrm{P}\big\{\mbox{Type I Error in $HT_j$}\big\} \leq \mathrm{P}\big\{\mbox{Reject } H_1^l\,|\,\pmb\mu \in H_1^l\big\}=\mathrm{P}\big\{\mbox{Type II Error in $HT_l$}\big\}\leq \alpha, \mbox{if } \pmb \mu\in H_1^l,
\end{align*}
or equivalently
\begin{align}\label{ErrorControl1}
\mathrm{P}\big\{\mbox{Type I Error in $HT_j$}\big\} \leq \alpha, \quad \forall\, \pmb \mu\in H_0^j,\, j=1,2,\dots,k.
\end{align}

Above all, we formulate the fixed-precision R\&S problem as a multiple HT problem in \eqref{test1} and illustrate that its precision (i.e., PCS guarantee in \eqref{PCS}) can be delivered by controlling the Type II error for each single $HT_j$ as presented in \eqref{ErrorControl}.

\subsubsection{The indifference-zone assumption}\label{subsec:IZ}
We next consider each $HT_j$ in \eqref{test1} individually, and notice that its Type I and II errors need to be controlled either directly or indirectly as discussed in Section~\ref{subsec:2general}. However, for a given set of samples, simultaneously controlling both types of error probabilities might be impossible. To show this, we connect these two error probabilities via the power function of the test, i.e.,
%\begin{small}
\begin{align*}
\beta_j(\pmb\mu) = \mathrm{P}\big\{\mbox{Reject } H_0^j\,|\,\pmb\mu\big\} = \left\{
\begin{aligned}
\mathrm{P}\big\{\mbox{Type I error in $HT_j$}\big\}, \quad &\mbox{ if } \mu_j\leq \max_{i\neq j}\mu_i, \\
1-\mathrm{P}\big\{\mbox{Type II error in $HT_j$}\big\}, \quad &\mbox{ if } \mu_j > \max_{i\neq j}\mu_i.
\end{aligned}
\right.
\end{align*}
%\end{small}
For most testing procedures, the power function $\beta_j(\pmb\mu)$ is continuous with respective to $\pmb\mu$. Then,
\begin{align}\label{conflict}
\mathrm{P}\big\{\mbox{Type I error in $HT_j$}\big\} = 1-\mathrm{P}\big\{\mbox{Type II error in $HT_j$}\big\}, \mbox{ when } \mu_j = \max_{i\neq j}\mu_i.
\end{align}
Obviously, this equation conflicts with the constraints stated in \eqref{ErrorControl} and \eqref{ErrorControl1}. Therefore, the testing procedure satisfying \eqref{ErrorControl} may not exist. It further reveals that in R\&S problems, we may not be able to select the best with the desired precision, when the means are sufficiently close to each other.

To overcome this obstacle, \cite{bechhofer1954single} introduced a so-called indifference-zone (IZ) parameter $\delta>0$, which refers to the smallest mean difference worth detecting. Given the IZ, the R\&S problems are modified to select the best alternative, when all the inferior alternatives are outside the IZ of the best. Accordingly, the PCS guarantee in \eqref{PCS} is rewritten as
\begin{align}\label{PCS_IZ}
\mbox{PCS-IZ}(\pmb\mu) = \mathrm{P}\big\{\mbox{Select the best alternative}\,|\,\pmb\mu\big\}\geq 1-\alpha, \ \forall\ \pmb\mu \in \Theta^\delta,
\end{align}
where $\Theta^{\delta}=\big\{\pmb\mu: \mu_{[k]}-\delta> \mu_{[k-1]}\big\}$ is called the IZ. Following the same logic in Section~\ref{subsec:2general}, this R\&S {problem} can be reformulated as a multiple HT problem, that is
\begin{align}\label{test_iz}
(HT_j^\delta)\quad \quad H_0^{j,\delta}: \mu_j+\delta \leq \max_{i\neq j}\mu_i, \quad \mbox{versus} \quad H_1^{j,\delta}:  \mu_j-\delta > \max_{i\neq j}\mu_i, \quad \forall\, j=1,2,\dots,k.
\end{align}
We remark here that, for any mean vector $\pmb\mu\in\Theta^\delta$ of interest, either $H_0^{j,\delta}$ or $H_1^{j,\delta}$ is true, which ensures the test above is well-defined.
%Furthermore, the corresponding power function turns into
%\begin{small}
%\begin{align*}
%\beta_j^\delta(\pmb\mu) = \mathrm{P}\big\{\mbox{Reject } H_0^j\,|\,\pmb\mu\big\} = \left\{
%\begin{aligned}
%\mathrm{P}\{\mbox{Type I error in $H_j^\delta$}\}, \quad &\mbox{ if } \mu_j+\delta \leq \max_{i\neq j}\mu_i, \\
%1-\mathrm{P}\big\{\mbox{Type II error in $H_j^\delta$}\}, \quad &\mbox{ if } \mu_j-\delta > \max_{i\neq j}\mu_i.
%\end{aligned}
%\right.
%\end{align*}
%\end{small}

Given the IZ parameter $\delta$, the corresponding power function is defined in two non-adjacent sets, i.e., $\{\pmb\mu: \mu_j+\delta \leq \max_{i\neq j}\mu_i\}$ and $\{\pmb\mu: \mu_j-\delta > \max_{i\neq j}\mu_i\}$. This frees us from facing the adjacent point, at which the Type I and II error probabilities cannot be controlled as desired because their sum is forced to be one. Therefore, in presence of the IZ parameter, it becomes possible to control both types of errors for each $HT_j^\delta$, or the Type II errors for the $HT_j^\delta (j=1,2,\dots,k)$. Accordingly, the R\&S problems with PCS-IZ can also be tackled. In Section~\ref{sec:fixed-precision}, we will explain in detail how several representative R\&S procedures are derived along this line.

\subsubsection{PCS and PGS}
As stated in Section~\ref{subsec:IZ}, the PCS guarantee in \eqref{PCS} is difficult to deliver, therefore the IZ parameter is introduced and the R\&S problems are restricted to a smaller mean vector space. As a consequence, the PCS-IZ guarantee in \eqref{PCS_IZ} is delivered whenever the best mean is at least $\delta$ larger than the others. However, in practice, several alternatives may have means that fall into the indifference zone, and these alternatives are called good alternatives. According to the definition of IZ, we should be indifferent if one of these good alternatives is selected as the best. Hence, we may care about the probability of good selection (PGS) rather than the original PCS, where the PGS guarantee is represented as
\begin{align}\label{PGS}
\mbox{PGS} {(\pmb\mu)}= \mathrm{P}\big\{\mbox{Select a good alternative}\,|\,\pmb\mu\big\}\geq 1-\alpha, \ \forall\ \pmb\mu \in \Theta.
\end{align}
{In the area of multi-armed bandits, a good selection is viewed as an approximately correct selection, and accordingly the PGS guarantee is also called the probably approximately correct (PAC) selection guarantee \citep{even2006action, ma2017an}.}

Notice that for the R\&S procedures with PCS-IZ guarantee, it is natural to expect that they could also deliver the PGS guarantee. Unfortunately, several counterexamples have been provided \citep{eckman2018guarantees}.

In the following, we attempt to explain this phenomenon from the hypothesis-testing perspective. Similar to Section~\ref{subsec:2general}, to select a good alternative, it suffices to test, for any given alternative $j$, whether it is a good alternative, i.e., $\mu_k+\delta>\max_{i\neq k} \mu_i$. Therefore, we formulate the R\&S problems with PGS guarantee as a multiple HT problem, that is
\begin{align}\label{test_gs}
(HT_j^G)\quad\quad \quad H_0^{j,G}: \mu_j +\delta\leq \max_{i\neq j} \mu_i \quad \mbox{ versus } \quad H_1^{j,G}: \mu_j +\delta >\max_{i\neq j} \mu_i,\ \forall\, j.
\end{align}

Suppose that a procedure with PCS-IZ guarantee of \eqref{PCS_IZ} exists, and we want to know whether it can deliver the PGS guarantee in \eqref{PGS}. According to the previous analysis, an easy way is by checking the Type II error constraints presented in \eqref{ErrorControl}. In Table~\ref{table1}, we summarize the R\&S problems with different probability guarantees and their corresponding HT formulations. Table~\ref{table1} shows that $H_0^{j,G} = H_0^{j,\delta}$. However, $H_1^{j,G}$ refers to a larger mean vector space than $H_0^{j,\delta}$. Therefore, the Type II error probability in $HT_j^{G}$ may not satisfy \eqref{ErrorControl} even though it is satisfied in $HT_j^\delta$. In other words, the PGS guarantee can not be guaranteed. {To overcome this drawback, \cite{eckmanPGS} constructed several sufficient conditions under which the PGS-IZ guarantee can imply the PGS guarantee.}

\begin{table}[h]
\centering
\caption{R\&S and their HT formulations.}\label{table1}
\begin{tabular}{ccccc}
\hline
Goal of R\&S && Means && HT Formulations\\
\hline
\multirow{2}{*}{PCS} && \multirow{2}{*}{$\mu_{[k]}>\mu_{[k-1]}$} && \multirow{2}{*}{$H_0^{j}: \mu_{j} \leq \max\limits_{i\neq j}\mu_i, \mbox{ v.s. } H_1^{j}:  \mu_{j} > \max\limits_{i\neq j}\mu_i, \  \forall\, j$}\\
&&&&\\
\multirow{2}{*}{PCS-IZ} && \multirow{2}{*}{$\mu_{[k]}-\delta>\mu_{[k-1]}$} && \multirow{2}{*}{$H_0^{j,\delta}: \mu_{j} +\delta\leq \max\limits_{i\neq j}\mu_i, \mbox{ v.s. } H_1^{j,\delta}:  \mu_{j} -\delta> \max\limits_{i\neq j}\mu_i, \  \forall\, j$}\\
%PCS && $\mu_{[k]}-\delta >\mu_{[k-1]}$ && $H_0^j: \mu_{j}+\delta \leq \max_{i\neq j}\mu_i, \mbox{ v.s. } H_1^j:  \mu_{j}-\delta > \max_{i\neq j}\mu_i, \forall\, j$\\
&&&&\\
\multirow{2}{*}{PGS} && \multirow{2}{*}{$\mu_{[k]}>\mu_{[k-1]}$} && \multirow{2}{*}{$H_0^{j,G}: \mu_{j} +\delta \leq \max\limits_{i\neq j}\mu_i, \mbox{ v.s. } H_1^{j,G}:  \mu_{j} +\delta> \max\limits_{i\neq j}\mu_i, \  \forall\, j$}\\
&&&&\\
\hline
\end{tabular}
\end{table}

On the opposite side, Table \ref{table1} depicts that PGS guarantee implies the PCS-IZ guarantee. Thus, interest has recently emerged in developing the procedures with PGS guarantee, such as \cite{fan2016indifference}, \cite{eckman2018guarantees}.

\subsection{Fixed-budget R\&S}\label{subsec:DP2}
In this section, we consider the R\&S procedures under a fixed sampling budget. By its nature, one can always select the alternative with the largest sample mean as the best when the sampling budget is exhausted. Therefore, the key issue here is how to allocate the budget efficiently. When the allocations can be made multiple times, one effective method is to re-determine the allocation adaptively at each stage based on the sampling information collected so far. Thus, a dynamic-programming \citep{bellman1966dynamic, bertsekas1995dynamic} formulation looks proper to derive an optimal allocation policy.

Under the DP formulation, R\&S problems turn into finding a sequence of sampling allocation decisions to optimize the precision of the final selection. Besides the PCS used in Section~\ref{sec:HT}, another popular measure to describe the precision of selection is expected opportunity cost (EOC). In fact, PCS is related to the so-called 0-1 loss, i.e., only a correct selection acquires a reward, while EOC describes the precision of selection by its opportunity cost. Particularly, when EOC is used, a non-best selection also gets a reward proportional to the discrepancy in the mean from the best one, which corresponds to a linear loss function. Instead of focusing on the final selection, some researchers have chosen to optimize the way how the  information has been collected, e.g., by maximizing the expected value of information (EVI) collected at each stage.

\subsubsection{Fixed-budget R\&S  formulated as dynamic programming}
Suppose that a total sampling budget $N$ is allocated to the $k$ alternatives progressively along $T$ stages, each endowed with a budget of $\tau=N/T$. (In the special case when $\tau=1$, the samples are allocated one by one.) {Assume that, at each stage $t$ ($t=1,2,\dots,T$), the $\tau$ samples are collected according to some sampling allocation policy, termed by $\pi_t$.} Apparently, the information about the alternatives is revealed gradually along the sequential sampling. To track the process, we denote $\mathcal E_0$ the initial information on the alternatives and $\mathcal E_t$ the information collected up to the end of stage $t$, for $t=1,2,\dots,T$. {The inter-stage updating rule of the information can be defined by a transition function $f_t$, i.e., $\mathcal E_{t}=f_t(\mathcal E_{t-1},\pi_t,\xi_t)$, where $\xi_t$ refers to the randomness of the samples collected at stage $t$.} After the final stage,  the selection decision is made based on all the information (i.e., $\mathcal E_T$) that is collected.

Let $V(\mathcal E_T)$ denote the terminal value function we want to optimize. For instance, when our objective is to minimize the probability of incorrect selection (i.e., $1-$PCS), the value function can be set as the 0-1 function which is 1 if the selected alternative is not the best and 0 otherwise. Then, the R\&S procedures are formulated as a DP, which is
\begin{equation}\label{dp}
  \min_\pi \mathrm E^\pi[V(\mathcal E_T)],
\end{equation}
where the decision is a sequence of {allocation} policies, i.e., $\pi=(\pi_1,\ldots,\pi_{T})$. In the literature, the DP problem is often handled recursively through the associated Bellman equation,
\begin{align}\label{eqn:bellman}
V^*_t(\mathcal E_t)=\min_{\pi_{t+1}}\mathrm E[V^*_{t+1}(\mathcal E_{t+1})], \quad t=T-1, T-2, \ldots, 0.
\end{align}
where the value function $V_t^*(\mathcal E_t)$ defines the optimal expected cost-to-go from current stage $t$ to the terminal and the terminal cost $V^*_T(\mathcal E_T)=V(\mathcal E_T)$.

Notice that the Bellman equation builds the relationship between the value functions in the current and next stages. As a consequence, the original DP is broken into a series of static optimization problems, although in a stage-by-stage and recursive form. However, in practice, the Bellman equation is  typically difficult to solve and the difficulty is illustrated as follows. To solve the Bellman equation, the next-to-terminal cost-to-go $V_{t+1}^*(\mathcal E_{t+1})$ in \eqref{eqn:bellman} has to be calculated by backward iterations. Unfortunately, these calculations tend to be increasingly difficult as the number of stages increases due to the ``curse of dimensionality".  In Section~\ref{sec:DP}, we will explain in detail how existing studies have resolved this problem and obtained the corresponding sample allocation rules (or R\&S procedures).

\subsubsection{Consistency of fixed-budget procedures}
With a fixed sampling budget, the DP R\&S procedures provide no probability guarantee on the correctness of selection. Alternatively, they usually process another appealing property of consistency. A procedure is said to be consistent if its selected alternative will converge to the true best as the total budget goes to infinity.

The consistency of a DP procedure is generally difficult to show directly. As long as all the alternatives receive infinite sampling budget in the limit, we will always have the exact information on the ranking of their true means to select the best correctly. Hence, asymptotically infinite samples on all the alternatives often works as a sufficient condition to verify the consistency of a procedure in the literature.

\subsection{Connection to the frequentist and Bayesian formulations}
Before this paper, the R\&S procedures under fixed-precision and fixed-budget were often classified into the frequentist and Bayesian procedures in the literature \citep{KimNelson2006b}. The main reason is that the precision of selected alternative or generally the value function in DP is often described under the corresponding frequentist or Bayesian probability models. However, there are some exceptions. For instance, \cite{frazier2014fully} proposed a R\&S procedure with PCS guarantee under a Bayes-inspired framework and \cite{chen2000simulation} suggested a R\&S procedure with a fixed budget under a frequentist framework.

Moreover, given that the R\&S problems under fixed-precision can be formulated as a hypothesis test, any testing rule, frequentist or Bayesian, can ideally be used to derive the corresponding R\&S procedures. Similarly, more sample allocation (or R\&S) procedures can be derived under either a frequentist or Bayesian framework for the R\&S problems under a fixed sampling budget. Therefore, in our view, R\&S procedures can be properly classified based on their underlying methodical formulations (i.e., HT or DP).

\section{Fixed-precision procedures}\label{sec:fixed-precision}
Considering the fixed-precision constraint, most of the existing R\&S procedures are designed under the IZ formulation and deliver the PCS-IZ guarantee in \eqref{PCS_IZ}.  These procedures are often called IZ procedures. Following the discussion in Section~\ref{sec:HT}, we will first show in detail how the stage-wise and sequential IZ procedures are derived by addressing the corresponding HT problem \eqref{test_iz}. Then we move to the newly designed IZ-free procedure which is able to deliver both the PCS and PGS guarantees.

%Notice that \eqref{test_iz} is not easy to handle, because this consists of $k$ simultaneous HTs with each dealing with the pairwise comparison between alternative $j$ and one of the other $k-1$ alternatives. To resolve this difficulty,
%%%\begin{small}
%%%\begin{eqnarray}\label{test2}
%%%\begin{split}
%%%H_0^1: \mu_{k}\leq \mu_{1} \quad & \mbox{versus}& \quad  H_1^1: \mu_{k}> \mu_{1}. \\
%%%H_0^2: \mu_{k}\leq \mu_{2} \quad &\mbox{versus}& \quad  H_1^2: \mu_{k}> \mu_{2}. \\
%%%& \vdots &\\
%%%H_0^{k-1}: \mu_{k}\leq \mu_{k-1} \quad& \mbox{versus}& \quad   H_1^{k-1}: \mu_{k}> \mu_{k-1}.
%%%\end{split}
%%%\end{eqnarray}
%%%\end{small}
%\begin{align}\label{test_iz_1}
%H_{0i}^\delta: \mu_{k}+\delta\leq \mu_{i} \quad \mbox{versus} \quad  H_{1i}^\delta: \mu_{k}-\delta> \mu_{i}, \quad \mbox{ for } i=1,2,\dots,k-1.
%\end{align}
%In doing so, the original multiple test is broken into a group of single tests. Then, we only need to handle these $k-1$ single tests separately, on which there is a vast volume of literature. Given the procedure for each single test, one can integrate them together to construct a procedure for the original \eqref{test_iz} as follows: reject $H_0^\delta$ if all the $H_{0i}^\delta\, (i=1,2,\dots,k-1)$ are rejected, and accept $H_0^\delta$ if at least one $H_{0i}^\delta$ is accepted.

Before moving to the next part, we first set up some notations. Let $X_{ij}$ denote the $j$th observation from alternative $i$, for $i=1,2,\dots,k$ and $j=1,2,\dots$. Unless specifically stated, we assume these observations are independent across alternatives and  $\{X_{ij}: j=1,2,\dots\}$ are independent and identically distributed (i.i.d.) Gaussian distribution with mean $\mu_i$ and variance $\sigma_i^2$. Let $\bar X_i(n)$ and $S_i^2(n)$ denote the sample mean and sample variance calculated based on the first $n$ samples from alternative $i$.

%
%Applying the well-known and commonly used Z-test procedure for the single hypothesis test in \eqref{test2}, the testing procedure for the original\eqref{test1} is constructed as follows,
%\begin{align*}
%\mbox{reject } H_0, &\quad\mbox{ if }\quad \bar{X}_k(n)-\bar{X}_i(n)\geq h\sqrt{2\sigma^2/n}, \ \forall\, i\neq k;\\
%\mbox{accept } H_0, &\quad\mbox{ if }\quad \bar{X}_k(n)-\bar{X}_i(n)\leq -h\sqrt{2\sigma^2/n},\  \exists\, i\neq k,
%\end{align*}
%where $n, \sigma^2$ is the number and variance of samples collected from each alternative, and  $h$ is  chosen to control the Type I error as shown in \eqref{ErrorControl}. However, such $h$ may not exist. Notice that
%\begin{align*}
%\mathrm{P}\{\mbox{Type I error}\} = \mathrm{P}\left\{\bar{X}_k(n)-\bar{X}_i(n)\geq h\sqrt{2\sigma^2/n}, \ \forall i\neq k\,\big|\, \mu_k< \max_{i\neq k} \mu_i\right\}
%\end{align*}
%When all the means are the same, i.e., $\mu_1=\dots=\mu_k$, it is intuitive that the statistical test cannot tell their differences and thus views each mean as the largest with the same probability $1/k$. In other words, it rejects a given mean as the best with  probability $(k-1)/k$ and this is the reason why the second equality above holds.

\subsection{Stage-wise R\&S procedures}\label{subsec:stagewise}
We start from deriving Bechhofer's procedure \citep{bechhofer1954single}, which is probably known as the first R\&S procedure in the literature. It considers a special case where the variances across all alternatives are common and known, i.e., $\sigma_1^2=\sigma_2^2=\dots=\sigma_k^2=\sigma^2$, and the goal is to deliver the PCS-IZ guarantee. In this case, one natural procedure for its corresponding HT problem in \eqref{test_iz} works as follows. For $j=1,2,\dots,k$,
\begin{align*}
\mbox{reject } H_0^{j, \delta} \mbox{ (i.e., select alternative $j$)}, &\quad\mbox{ if }\quad \bar{X}_j(n)-\max_{i\neq j}\bar{X}_i(n)\geq z,
\end{align*}
and accept $H_0^{j, \delta}$ otherwise. Here the constant $z$ and the common sample size $n$ of all alternatives need to be carefully chosen.

Only a single alternative is expected to be returned as the best. Straightforwardly, it occurs if only one $H_0^{j,\delta}$ is rejected. This suffices to require that the rejection regions for $H_0^{j,\delta}\, (j=1,2,\dots,k)$ compose the disjoint partition of the whole space $\real_+^k$. One way to achieve this goal is setting $z=0$. In doing so, the alternative with the largest sample mean is selected as the best. Besides, the common sample size $n$ is chosen such that the Type II error probability for each $HT_j^\delta$ satisfies \eqref{ErrorControl}, and specifically,
\begin{align}\label{Bechhofer_typeI}
\notag &\quad\; \mathrm{P}\big\{\mbox{Type II error in $HT_j^{\delta}$}\big\} \\
\notag &=\mathrm{P}\left\{\bar{X}_j(n)-\max_{i\neq j}\bar{X}_i(n)< 0\,\bigg|H_1^{j,\delta}\right\}\\
\notag & = \mathrm{P}\left\{\max_{i\neq j}\frac{\sqrt{n}[\bar X_i(n)-\bar X_j(n)-(\mu_i-\mu_j)]}{\sqrt{2\sigma^2}}> -\max_{i\neq j}(\mu_i-\mu_j)\sqrt{\frac{n}{2\sigma^2}}, \forall i\neq j\,\bigg|H_1^{j,\delta}\right\}\\
& \leq \mathrm{P}\left\{\max_{i\neq j} Z_i > \delta \sqrt{\frac{n}{2\sigma^2}}\right\}\leq \alpha,
\end{align}
%\begin{align}\label{Bechhofer_typeI}
%\notag \mathrm{P}\big\{\mbox{Type II error in $HT_j^{\delta}$}\big\} &=\mathrm{P}\left\{\bar{X}_j(n)-\max_{i\neq j}\bar{X}_i(n)< 0\,\bigg|H_1^{j,\delta}\right\}\\
%\notag & \leq \mathrm{P}\left\{\max_{i\neq j}\frac{\sqrt{n}[\bar X_i(n)-\bar X_j(n)-(\mu_i-\mu_j)]}{\sqrt{2\sigma^2}}> \delta\sqrt{\frac{n}{2\sigma^2}}, \forall i\neq j\,\bigg|H_1^{j,\delta}\right\}\\
%& = \mathrm{P}\left\{\max_{i\neq j} Z_i > \delta \sqrt{\frac{n}{2\sigma^2}}\right\}\leq \alpha,
%\end{align}
where $(Z_i,\, i\neq j)$ is a $(k-1)$-dimensional multivariate Gaussian random variable with means 0, variances 1 and common pairwise correlations 1/2.  Let $h$ denote the $(1-\alpha)$ quantile of the maximum of $Z_i\, (i\neq j)$. The common sample size $n$ is chosen as
\begin{align}\label{bechhofer_n}
n = \left\lceil\frac{2h^2\sigma^2}{\delta^2}\right\rceil,
\end{align}
where $\lceil x \rceil$ denotes the smallest number no smaller than $x$.

Following the testing procedure above, a R\&S procedure can be constructed. It first determines the common sample size allocated to each alternative as \eqref{bechhofer_n}. Then it selects the alternative with the largest sample mean as the best. This is exactly Bechhofer's procedure.

Regarding Bechhofer's procedure, we make two remarks here.
\begin{itemize}
\item[(i)] From \eqref{Bechhofer_typeI}, we see that the worst-case of Type II error probabilities is attained when the best mean is exactly $\delta$ better than all the others, i.e., $\mu_{[k]}-\delta = \mu_{[k-1]}=\dots=\mu_{[1]}$. Thus, this configuration of means is the most difficult situation in $\Theta^\delta$ and \cite{bechhofer1954single} names it the least favorable configuration (LFC) of means.
\item[(ii)] Bechhofer's procedure is also able to deliver the PGS guarantee in \eqref{PGS}. To verify this statement, we only need to prove that the Type II error constraint in \eqref{ErrorControl} can be achieved while applying the procedure to address the $HT_j^G$ for all $j$. This proof is easily accomplished and therefore omitted in this study.
\end{itemize}

\cite{rinott1978two} extends Bechhofer's procedures to the situation where the variances across alternatives are unknown and unequal. To handle this situation, Bechhofer's procedure is modified in three aspects. First, an initial stage is included in which a small number of samples are generated to estimate the unknown variances. Second, the total sample sizes allocated to each alternative are not the same any more but are set to be positively proportional to its sample variance. Third, the constant $h_R$ in the total sample size $N_i$ needs to be modified accordingly. {Finding this constant needs to solve a root-find problem with integration, i.e.,
\begin{align*}
\int_{-\infty}^\infty \Psi^{k-1}_{n_0-1}(t+h_R)\psi_{n_0-1}(t) dt = (1-\alpha),
\end{align*}
where $\Psi_{n_0-1}$ and $\psi_{n_0-1}$ denote the cumulative distribution function and probability density function of a standard student-t distribution with $n_0-1$ degrees of freedom, respectively. } Historically, due to the limited computational capacity, it is considered difficult to solve; hence, tables are provided \citep{wilcox1984table, santner1995design, goldsman1998comparing}. The new two-stage procedure (named as the Rinott's procedure) is presented as follows.

\begin{algorithm}
\caption{Rinott's Procedure}
\vspace{3pt}
\textbf{Require:} Number of alternatives $k$, common first-stage sample size $n_0\geq 2$, PCS $1-\alpha$, IZ parameter $\delta$, a constant $h_R$.
\begin{algorithmic}[1]
    \State Generate $n_0$ samples for each alternative $i$ and calculate the sample variance $S_i^2(n_0)$.
    \vspace{3pt}
    \For{$i\leftarrow 1:n$}
    \vspace{3pt}
        \State Let
        \begin{small}
        \begin{align}\label{Rinott_N}
        N_i\leftarrow \max\left\{n_0,\left\lceil\frac{h_R^2 S_i^2(n_0)}{\delta^2} \right\rceil\right\}.
        \end{align}
        \end{small}
        \State Generate $N_i-n_0$ samples from alternative $i$ and calculate the sample mean $\bar{X}_i(N_i)$.
    \EndFor
    \vspace{3pt}
    \State Select $\arg\max_{i=1,2,\dots,k}\bar{X}_i(N_i)$ as the best.

\end{algorithmic}
\end{algorithm}

As the simplest and most popular IZ procedure, there are a lot of variations of Rinott's procedure. For instance, to avoid the complexity in calculating $h_R$, some procedures \citep{clark1986bonferroni} adopt Bonferroni's inequality and set it approximately as the $1-\alpha/(k-1)$ quantile of a t-distribution with $n_0-1$ degrees of freedom \citep{banerjee1961confidence}. As a price, it often leads to more conservativeness, which means that a larger sample size is needed for the procedure. Another variation of Rinott's procedure worth mentioning is the use of common random numbers (CRNs) \citep{clark1986bonferroni, nelson1995using}. CRNs artificially introduce a positive correlation between the observations from each pair of alternatives, thus {decreasing} the variance of their sample mean difference. In doing so, the R\&S process becomes much easier and the sample size required is ultimately reduced.

\subsection{Sequential R\&S procedures}\label{subsec:sequential}
Paulson's procedure is one of the early sequential R\&S procedures, and this subsection will start from re-deriving this procedure from the hypothesis-testing perspective. Same as Bechhofer's procedure, Paulson's procedure also considers the special case with common and known variances, i.e., $\sigma_1^2 =\sigma_2^2=\dots=\sigma_k^2=\sigma^2$.

Similar to Section~\ref{subsec:stagewise}, we first consider each $HT_j^\delta$ individually and our task is to design a sequential testing procedure for it. However, such sequential procedure is not trivial because it involves multiple pairwise comparisons between alternatives. As a remedy, we break down $HT_j^\delta$ into a group of HT problems, each of which considers a pairwise comparison between alternative $j$ and one of the other alternatives. Particularly, $HT_j^\delta$ is decomposed into
\begin{align}\label{decomposition}
(HT_{ji}^\delta) \quad\quad H_0^{ji,\delta}: \mu_j+\delta\leq \mu_i\quad \mbox{ versus }\quad H_1^{ji,\delta}: \mu_j-\delta>\mu_i, \forall i\,\neq j.
\end{align}
Meanwhile, to control the Type II error in  $HT_j^\delta$ at most $\alpha$ as desired in \eqref{ErrorControl}, we adopt Bonferroni's inequality and require
\begin{align}\label{PairwiseError}
\mathrm{P}\big\{\mbox{Type II error in $HT_{ji}^{\delta}$}\big\} \leq \alpha/(k-1), \forall\, i\neq j.
\end{align}
A sequential procedure for $HT_{ji}^\delta$ is noticeably easy to obtain while satisfying \eqref{PairwiseError}, and a vast volume of literature supports it. Specifically, we may use Wald's sequential probability ratio test (SPRT) \citep{wald1945sequential, wald2004sequential} which,
\begin{align*}
\mbox{rejects } H_0^{ji, \delta}, \quad &\mbox{if}\quad n(\bar{X}_j(n)-\bar{X}_i(n))\geq a-\lambda n,\\
\notag \mbox{accepts } H_0^{ji,\delta}, \quad &\mbox{if}\quad n(\bar{X}_j(n)-\bar{X}_i(n))\leq -a+\lambda n,
\end{align*}
and continues to take samples otherwise. Here $0<\lambda <\delta$ and $a$ is chosen as $a = \ln\left(\frac{k-1}{\alpha}\right)\frac{\sigma^2}{\delta-\lambda}$.

Now the original R\&S problem is reformulated as $k(k-1)$ simultaneous HT problems, i.e., $HT_{ji}^\delta$, for $j\neq i$. Each $HT_{ji}^\delta$ considers that the pairwise comparison between alternatives $j$ and $i$ and is resolved by a sequential procedure as mentioned above. Intuitively, at any time of the sampling process, we should select alternative $j$ as the best if all the $H_0^{ji, \delta}\, (i\neq j)$ are rejected. We eliminate alternative $j$ from consideration if one of the $H_0^{ji, \delta}\, (i\neq j)$ is accepted. Otherwise we continue to take samples otherwise. Once an alternative  is eliminated, we should stop taking samples from this alternative and {abandon} all the $HT_{ji}^\delta$ regarding it. For clarity, $I(n)$ denotes the set of surviving alternative right before stage $n$, then a sequential procedure is designed as
\begin{align*}
\mbox{selecting alternative $j$}, \quad &\mbox{if}\quad n(\bar{X}_j(n)-\bar{X}_i(n))\geq a-\lambda n, \forall\,i\in I(n)\mbox{ and }i\neq j,\\
\notag \mbox{eliminating alternative $j$ }, \quad &\mbox{if}\quad n(\bar{X}_j(n)-\bar{X}_i(n))\leq -a+\lambda n, \exists\, i\in I(n)\mbox{ and }i\neq j.
\end{align*}
It continues to take samples on the surviving alternatives otherwise. This sequential procedure is known as Paulson's procedure.

\cite{kim2001fully} extend Paulson's procedure to the case of unknown and unequal variances. Similar to the previous two-stage procedures, Kim and Nelson’s ($\mathcal{KN}$) procedure also uses an additional initial stage of sampling to estimate the unknown variances. After the variances are estimated, it then starts screening alternatives just as Paulson's procedure does. In addition, replacing Paulson's bound by a tighter bound of \cite{Fabian1974} and considering the estimated variances which are random variables, $\mathcal{KN}$ procedure re-assigns the values of $\lambda$ and $a$ to ensure the same PCS guarantee. The detailed $\mathcal{KN}$ procedure is presented in Procedure~\ref{proc:KN}.

\begin{algorithm}[htp]
\caption{$\mathcal{KN}$ Procedure}\label{proc:KN}
\vspace{3pt}
\textbf{Require:} Number of alternatives $k$, common first-stage sample size $n_0\geq 2$, PCS $1-\alpha$, IZ parameter $\delta$, a constant $h$.
\vspace{3pt}
\begin{algorithmic}[1]
    \State Set
    \begin{footnotesize}
    \[
    \eta = \frac{1}{2}\left[\left(\frac{2\alpha}{k-1}\right)^{-2/(n_0-1)}-1\right].
    \]
    \end{footnotesize}
    \State $I\leftarrow \{1,2,\dots,k\}$, $h^2 = 2\eta(n_0-1)$, $n\leftarrow n_0$.
    \State Generate $n_0$ samples to each alternative $j$ and calculate $\bar{X}_i(n_0)$. For $i, j\in I$,
    \begin{footnotesize}
    \[
    S_{ji}^2 = \frac{1}{n_0-1}\sum_{l=1}^{n_0} \left[X_{jl}-X_{il}-(\bar{X}_j(n_0)-\bar{X}_i(n_0))\right]^2.
    \]
    \end{footnotesize}
    \While{$|I|>1$}
        \State Set $\small{W_{ji}=\max\Big\{0,\frac{\delta}{2n}\Big(\frac{h^2S_{ji}^2}{\delta^2}-n\Big)\Big\}}$ and
        \begin{small}
        \[
        I\leftarrow \big\{j: j\in I \mbox{ and } \bar{X}_j(n)-\bar{X}_i(n)\geq -W_{ji}(n), \forall i\in I, i\neq j\big\}
        \]
        \end{small}
        \State Take an additional observation from each alternative $j\in I$, and set $n\leftarrow n+1$.
    \EndWhile
    \State Select the alternative in $I$ as the best.
\end{algorithmic}
\end{algorithm}

\begin{figure}[ht]
\centering
\includegraphics[width=0.65\textwidth]{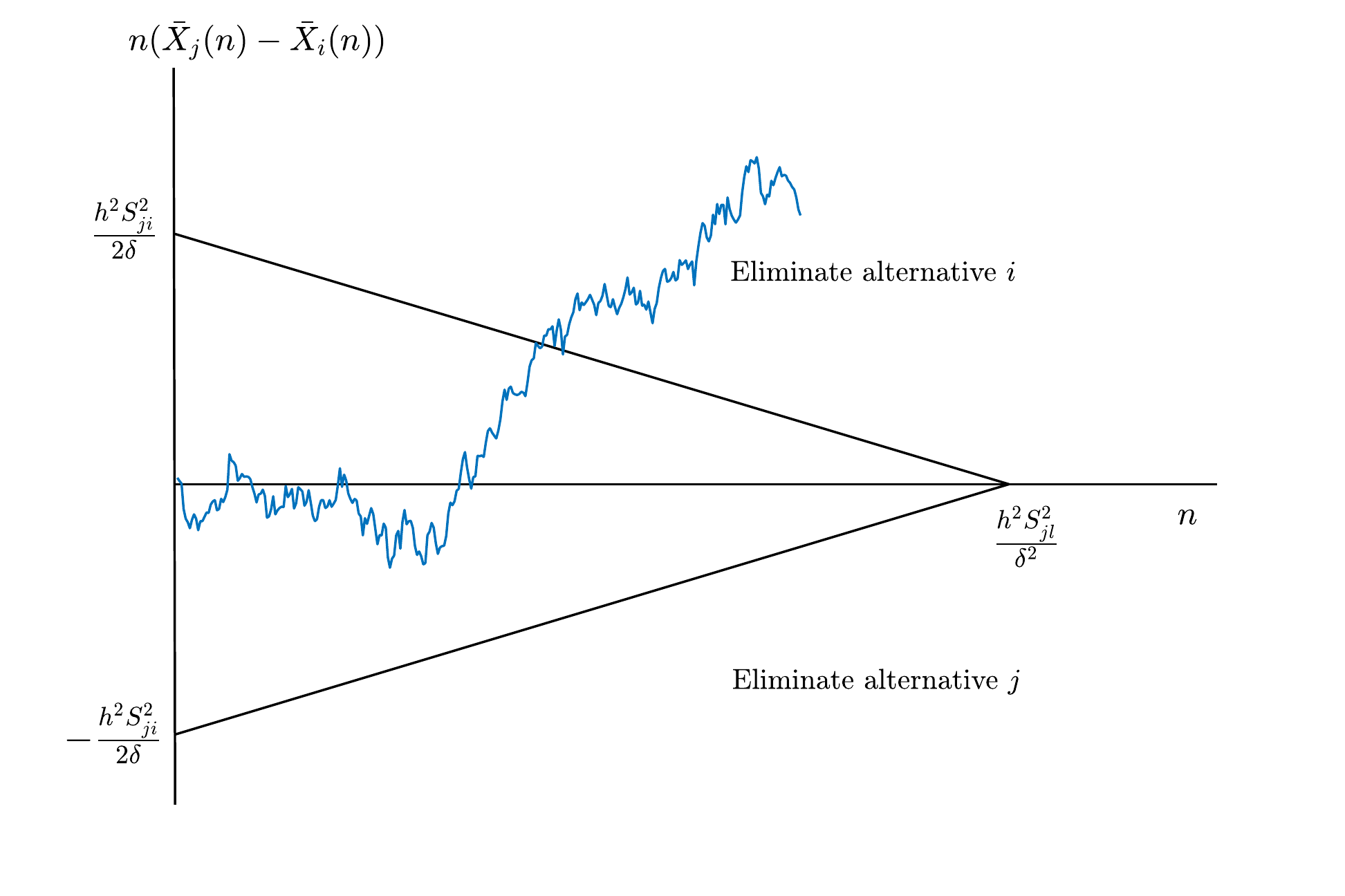}
\caption{Triangular Region for the $\mathcal{KN}$ procedures.}\label{fig:kn}
\end{figure}

An intuitive way to understand $\mathcal{KN}$ procedure is presented in Figure \ref{fig:kn}. For each pair of alternatives $j$ and $i$, it constructs the partial-sum process of their mean difference $\{n(\bar{X}_j(n)-\bar{X}_i(n)): n=1,2,\dots\}$. Then, at each stage $n$, $\mathcal{KN}$ checks whether this partial-sum process exits from the triangular region and makes decisions accordingly .

The $\mathcal{KN}$ procedure has numerous variations, and this family of procedures is shown to be effective among IZ procedures \citep{KimNelson2006b,branke2007selecting}. All these variations are classified into two categories. The first category intends to enhance the efficiency of the $\mathcal{KN}$ procedure. For instance, \cite{Hong2006} designed a variance-dependent sampling rule. Moreover, \cite{tsai2009fully} and \cite{TsaiLuo2017} adopted the control-variates technique. In another study, \cite{NSGS2001} took advantage of the first-stage samples to screen out alternatives that are unlikely to be the best. The second category intends to address different practical situations. For instance, \cite{hong2005tradeoff} considered the cost of switching between alternatives to take samples and designed a new procedure to balance the tradeoff between sampling and switching costs. In a follow-up study, \cite{hong2007selecting} noticed a situation where alternatives may be revealed sequentially, thus designing a new procedure for this situation. Meanwhile, \cite{kim2006asymptotic} studied the steady-state experiments and designed a new procedure achieving the PCS guarantee asymptotically.

\subsection{Indifference-zone-free R\&S procedures}
In the previous Sections~\ref{subsec:stagewise} and \ref{subsec:sequential}, we have seen how the IZ formulation (i.e., $\pmb\mu\in\Theta^\delta=\{\pmb\mu: \mu_{[k]}-\delta>\mu_{[k-1]}\}$) helps to achieve the PCS guarantee. However, problem remains whether a R\&S procedure with the PCS guarantee can be developed for all possible mean vectors in $\Theta$.

To solve this problem, \cite{fan2016indifference} proposed an IZ-free procedure. We call it FHN procedure and present it as Procedure 3. Similar to $\mathcal{KN}$ procedure, it decomposes a R\&S problem into a group of pairwise comparisons and designs a procedure for each pairwise comparison. When $\pmb\mu\in \Theta$, the pairwise mean differences might be arbitrarily close to zero. Then, the desired procedure is intended to detect whether these mean differences are zero or not. Motivated by the law of iterated logarithm, this IZ-free procedure adopts a new continuation region whose boundary function grows to infinity at a rate between $O(\sqrt{n\log\log n})$ and $O(n)$. For instance, a boundary function $\sqrt{[c+\log(n+1)](n+1)}$ is used as shown in Procedure 3.

\vspace{6pt}
\begin{algorithm}[htp]
\caption{FHN Procedure}
\vspace{3pt}
\textbf{Require:} Number of alternatives $k$, common first-stage sample size $n_0\geq 2$, PCS $1-\alpha$.
\vspace{3pt}
\begin{algorithmic}[1]
    \State Set
    \begin{footnotesize}
    \[
    c = -2\log\left(\frac{2\alpha}{k-1}\right).
    \]
    \end{footnotesize}
    \State $I\leftarrow \{1,2,\dots,k\}, n\leftarrow n_0$.
    \State Generate $n_0$ samples to each alternative $j$ and calculate $\bar{X}_j(n_0)$. For $i, j\in I$,
    \begin{footnotesize}
    \[
    S_{ji}^2(n_0) = \frac{1}{n_0-1}\sum_{l=1}^{n_0} \left[X_{jl}-X_{il}-(\bar{X}_j(n_0)-\bar{X}_i(n_0))\right]^2.
    \]
    \end{footnotesize}
    \While{$|I|>1$}
    \vspace{3pt}
        \State Set $\small{t_{ji}(n)=n/S_{ji}^2(n)}$ and $\small{g_{ji}(t_{ji}(n))=\sqrt{[c+\log(t_{ji}(n)+1)](t_{ji}(n)+1)}}$, and let
        \begin{small}
        \[
        I\leftarrow \{j: j\in I \mbox{ and } t_{ji}(n)[\bar{X}_j(n)-\bar{X}_i(n)]\geq -g_{ji}(t_{ji}(n)), \forall i\in I, i\neq j\},
        \]
        \end{small}
        \State Take an additional observation from each alternative $j\in I$, and set $n\leftarrow n+1$.
    \EndWhile
    \State Select the alternative in $I$ as the best.
\end{algorithmic}
\end{algorithm}

Now we illustrate from the HT perspective why this IZ-free procedure is able to achieve the PCS guarantee in \eqref{PCS}. As mentioned in Section~\ref{subsec:IZ}, the challenge for the conventional IZ procedures is how to control the Type I and Type II errors in each $HT_j$ simultaneously when the second-best mean is arbitrarily close to the best. Specifically, \eqref{conflict} shows that we might lose such control at the point $\pmb\mu^0$ with $\mu_j^0=\max_{i\neq j}\mu_i^0$, which is caused by the continuity of the power function. The FHN procedure resolves this challenge by forcing its power function $\beta_j(\cdot)$ to be discontinuous at $\pmb\mu^0$.

The FHN procedure addresses the $HT_j\, (j=1,2,\dots,k)$ by
\begin{align*}
\mbox{rejecting } H_0^{j}, \quad &\mbox{if }t_{ji}(n)[\bar{X}_j(n)-\bar{X}_i(n)]\geq g(t_{ji}(n)), \forall i\in I(n) \mbox{ and } i\neq j,\\
\notag \mbox{accepting  $H_0^{j}$}, \quad &\mbox{if }t_{ji}(n)[\bar{X}_j(n)-\bar{X}_i(n)]\leq -g(t_{ji}(n)), \exists\, i\in I(n) \mbox{ and } i\neq j,
\end{align*}
and continues sampling otherwise. Here $I(n)$ denotes the set of surviving alternatives right before stage $n$. Then, a careful derivation yields that
\begin{align*}
\beta_j(\pmb\mu)\geq 1-\alpha, \mbox{ for } \pmb\mu \mbox{ with }\mu_j>\max_{i\neq j}\mu_i, \mbox{ and }
\beta_j(\pmb\mu)\leq \alpha, \mbox{ for } \pmb\mu \mbox{ with }\mu_j\leq \max_{i\neq j}\mu_i,
\end{align*}
thereby demonstrating a discontinuous power function $\beta_j(\pmb\mu)$. The inequalities above also show that FHN procedure {satisfies} the constraints of error probability in \eqref{ErrorControl} and \eqref{ErrorControl1}, implying the desired PCS guarantee in \eqref{PCS} can be achieved.

\cite{fan2016indifference} also extended FHN procedure to incorporate an IZ parameter when it is available. Particularly, a stopping condition, based on the IZ parameter, is embedded into the original FHN procedure. The new procedure is shown to be able to achieve not only the PCS guarantee as in \eqref{PCS}, but also the PGS guarantee in \eqref{PGS}.

%In the previous Sections 3.1 and 3.2, a series of IZ procedures are introduced and they are able to achieve the PCS guarantee whenever the mean vector $\pmb\mu\in\Theta^\delta=\{\pmb\mu: \mu_{[k]}-\delta>\mu_{[k-1]}\}$. However, the original goal of PCS are required for all possible mean vector satisfying $\mu_{[k]}>\mu_{[k-1]}$.

% Conventionally, these procedures depend on the IZ parameter $\delta$.  One interesting question arising here is whether we can design an IZ procedure that is free of $\delta$. In doing so, this procedure is supposed to

\section{Fixed-budget procedures}\label{sec:DP}
In this section, we review the existing fixed-budget R\&S procedures related to the DP formulation. With a fixed sampling budget, the main task of R\&S procedures is to determine a sample allocation policy, which is formulated as a DP problem in \eqref{dp} as introduced in Section~\ref{subsec:DP2}. This DP problem is essentially a finite-horizon stochastic DP, and can be solved exactly by backward induction through Bellman equation \eqref{eqn:bellman}. However, this exact procedure is often impossible to execute due to the curse of dimensionality. This motivates the researchers to consider the suboptimal solutions generated by easily implementable approximation procedures. In particular, all the procedures reviewed in this section can be regarded as approximate dynamic programming (ADP) procedures.

\subsection{Static-allocation based procedures}
%\begin{itemize}
%\item Asymptotically optimal sampling ratio: Glynn and Juneja
%\item OCBA procedures
%\end{itemize}

%\subsubsection{Static Allocation and the Asymptotically Optimal Sampling Ratio}
As {a practically acceptable DP procedure is impossible to obtain}, one possible approach would be developing a good heuristic procedure instead. Intuitively, a superior DP procedure ``optimizes'' the way of collecting information about the mean of each alternative. Hypothetically, if we have perfect information at the beginning but still have to make selection based on the samples, a simple static allocation policy that maximizes the precision of the selection will be proper. For example, assuming the precision of selection is measured by the PCS guarantee in \eqref{PCS}, the optimal allocation policy can be determined by solving the following static optimization problem,
\begin{equation}\label{static}
%\max_{n_{[1]}+\cdots+n_{[k]}=N}\mathrm{PCS}\equiv
\max_{n_{[1]}+\cdots+n_{[k]}=N}\mathrm P\left(\bar X_{[k]}(n_{[k]}) > \max_{[j]\neq [k]}\bar X_{[j]}(n_{[j]})\right),
\end{equation}
where $n_{[i]}$ denotes the sample size allocated to alternative $[i]$, for $i=1,2,\dots,k$.
%This static problem seems much easier than the original DP problem.

Based on the static allocation policy, several procedures have been developed. The optimal computing budget allocation (OCBA) procedure, initiated by \cite{chen1996lower} and \cite{chen2000simulation}, is among the most famous static-allocation based procedures. Moreover, the OCBA procedure has also been extended to sequential settings, and the basic idea is to dynamically approximate the static allocation policy based on the sample information.

Taking the sequential algorithm of OCBA proposed by \cite{chen2000simulation} as an example, a total budget of $N$ is allocated to $T$ stages sequentially with each stage endowed with $\tau=N/T$. Perfect information is assumed in developing the OCBA procedure at first. Particularly, it assumes the information given at stage $t$ as $\mathcal E_t = \{(\mu_j,\sigma_j^2), j=1,\cdots k\}$ for $0\leq t\leq T$.
For any intermediate stage $t$, the allocation policy is determined by a static allocation problem as \eqref{static}, in which the budget for the first $t$ stages are reallocated for a myopic objective of maximizing PCS as if the selection is made at the end of the current stage.
\begin{equation*}
V_t^\text{OCBA}(\mathcal E_t)=\max_{n_{[1],t}+\cdots+n_{[k],t}=\tau t} \mathrm P\left(\bar X_{[k]}(n_{[k],t}) > \max_{ [j]\neq [k]}\bar X_{[j]}(n_{[j],t})\right).
\end{equation*}
Here $n_{[i],t}$ is the total sample size that is allocated to alternative $[i]$ up to the end of stage $t$, for $i=1,2,\dots,k$ and $t=1,2,\dots,T$. The allocation rule is then derived by approximating the PCS with Bonferroni's inequality and letting the budget per stage goes to infinity. The resulting allocation rule is presented in Step 5 in Procedure 4.

Moreover, using the large deviation theory, \cite{glynn2004large} derived the asymptotic optimal allocation policy for \eqref{static} that maximizes the exponential decay rate of the probability of incorrect selection as $N\rightarrow \infty$. Specially, they showed that the optimal allocation satisfies
\begin{align}\label{static_ratio}
\frac{n_{[i]}^*}{n_{[j]}^*}\approx\frac{\sigma_{[i]}^2/(\mu_{[k]}-\mu_{[i]})^2}{\sigma_{[j]}^2/(\mu_{[k]}-\mu_{[j]})^2},\quad \text{for }[i]\neq [j]\neq [k], \quad\mbox{ and }
n_{[k]}^*=\sigma_{[k]}\sqrt{\sum_{{[j]}\neq [k]}\left(\frac{n_{[j]}^*}{\sigma_{[j]}}\right)^2}.
\end{align}
This equation provides a theoretical benchmark on the optimality of static allocation policy. Careful investigation reveals that the optimal allocation coincides with the one in the OCBA procedure. Thus, the OCBA policy is asymptotically efficient.

%Now we are ready to derive several representative R\&S procedures via the approximated Bellman equation \eqref{eqn:bellman_myopic}. The first family of procedures is called the optimal computing budget allocation (OCBA) procedure \citep{chen1996lower, chen2000simulation}. \cite{chen2000simulation} considers the terminal value function as the 0-1 loss function,
%in attempt to maximize the PCS (or equivalently minimize the probability of incorrect selection). Accordingly, \eqref{eqn:bellman_myopic} is rewritten as
%\begin{align}
%1-\mbox{PCS} = \min_{n_{1,t}+n_{2,t}+\dots+n_{k,t}=b_t}\mathrm P\left\{U_{[k]}<U_{j}, \exists j\neq [k]\big|\mathcal E_t\right\},
%\end{align}
%where $[k]=\arg\max \mu_i$ denotes the index of the true best alternative and $b_t$ denotes total the sampling budget up to stage $t$.

%For the VFA, it first constructs an upper bound $\sum\limits_{j\neq [k]}\mathrm P\left\{U_{[k]}<U_{j}\right\}$ for the expected value function by the Bonferroni inequality, and then uses a normal approximation to describe the distribution of $(U_{[k]}-U_{j})|\mathcal{E}_t$ \citep{degroot2005optimal}. For the PSA, it optimizes the VFA obtained in an asymptotic regime as the sampling budget $\tau$ per stage grows to the infinity, and then derives the optimal sampling allocation policy $\{n_{1,t},n_{2,t},\dots, n_{k,t}\}$ as Step 3 in Procedure 3.

%This optimal policy is shown to be asymptotically efficient, in the sense that it can maximize the exponential decay rate of probability of incorrect selection \citep{glynn2004large}.

\begin{algorithm}
\caption{OCBA Procedure}
\vspace{3pt}
\textbf{Require:} Number of alternatives $k$, common first-stage sample size $n_0\geq 5$, total sampling budget $N$, sampling budget $\tau$ per stage.
\vspace{3pt}
\begin{algorithmic}[1]
    \State Generate $n_0$ samples from each alternative $i$.
    \State Set $t\leftarrow 0$, $n_{i,t}\leftarrow n_0$, $b_t\leftarrow\sum_{i=1}^k n_{i,t}$.
    \While{$b_t<N$}
    \vspace{3pt}
        \State Update {the sample mean} $\bar{x}_i$ \footnotemark and {the sample variance} $\hat\sigma_i^2$; $(k)\leftarrow \arg\max_i \bar{x}_i$ and $d_{(i)(k)}\leftarrow\bar{x}_{(k)}-\bar{x}_i$.
        \vspace{3pt}
        \State Set $b_{t+1}\leftarrow b_t+\tau$. Calculate the new budget allocation $n_{1,t+1}, n_{2,t+1},\dots, n_{k,t+1}$ satisfying $\sum_i n_{i,t+1}=b_{t+1}$ according to
        \begin{small}
        \begin{align*}
        \frac{n_{i,t+1}}{n_{j,t+1}}= \left(\frac{\hat\sigma_i/d_{(i)(k)}}{\hat\sigma_j/d_{(j)(k)}}\right)^2, \mbox{ for } i\neq j\neq (k), \mbox{ and } n_{(k),t+1} =\hat\sigma_{(k)}\sqrt{\sum_{i\neq (k)}n_{i,t+1}^2/\hat\sigma_i^2}.
        \end{align*}
        \end{small}
        \State Generate $ \max\{0, n_{i,t+1}-n_{i,t}\}$ samples from each alternative $i$. Set $t\leftarrow t+1$.
    \EndWhile
    \State Select $\arg\max \bar x_{i}$ as the best.
\end{algorithmic}
\end{algorithm}
\footnotetext{{Most of DP procedures are described as Bayesian procedures. Therefore, they are used to represent random variable with upper-case letters and represent their observations with lower-case letters. To keep in line with the existing literature, we use $\bar x$ in Section 4 to denote the observation of sample mean.}}

In practice we do not have perfect information about the means and variances of the alternatives, and the OCBA procedure suggests to use sample estimates instead, based on the available data at the beginning of each stage (see Step 4).
%It is worthy mentioning here that the OCBA procedure discussed above assumes that the samples among alternatives are independent and their variances are known. For the unknown variances in practice, it has to estimate the sample variances at the beginning of each stage .

Some variations of the above OCBA procedure have been proposed. \cite{he2007opportunity} adopted the linear loss function to measure the quality of the selection and designed an OCBA-type procedure; \cite{gao2017new} also consider the case of linear loss function but design an OCBA-type procedure based on the large-deviation theory. \cite{branke2007selecting} addressed the issue of unknown variances and proposed to use a student-t approximation; \cite{PengCCF2017} directly approximate the objective function in (4.17) using some feature functions. Moreover, \cite{fu2007simulation} considered the case of correlated samples across alternatives and showed that the optimal policy agrees with that of the independent case as the correlation vanishes.

\subsection{Two-stage approximation based procedures}
%\begin{itemize}
%\item OCBA assumes that the means and variances are known.
%\item The two-stage procedures take a first-stage sampling to estimate the unknown mean and variances, which are further used to construct the student approximation for the pairwise difference in the second stage.
%\item EI procedures (Chick)
%\end{itemize}

%Static-allocation based procedures like OCBA are developed by assuming the means and variances are known, and use the corresponding sample estimators to replace the unknown parameters in practical implementations. In this sense, the two-stage procedures are appropriate, which often include a first stage sampling to collect some information about the unknown parameters, and use the information to guide the second stage allocation decision. Based on the two-stage procedure, a sequential procedure can be immediately obtained by repeatedly executing the two-stage procedure. In the perspective of DP formulation, this kind of two-stage based procedure approximates the Bellman equation by myopically allocating budget to optimize the quality of selection at the current stage.

Static-allocation-based procedures like OCBA are developed by assuming that the means and variances are known. In addition, these procedures use the corresponding sample estimators to replace the unknown parameters in practical implementations. In contrast, another stream of research takes account of the unknown means and variances in developing the procedures. These procedures often contain two stages, which include a first-stage sampling to collect some information about the unknown parameters and then use the information to guide the second-stage allocation decision.

% In the perspective of DP formulation, this kind of two-stage based procedure approximates the Bellman equation by myopically allocating budget to optimize the quality of selection at the current stage.

As a representative, we shall review one famous two-stage procedure proposed by
\cite{chick2001new}, known as the expected value of information (EVI) procedure. In particular, we consider the one with linear loss and a budget constraint, or namely Procedure $\mathcal{LL(B)}$.
The procedure adopts Bayesian approach for updating the information collected about the mean performance of any alternative~$i$, which is assumed to be a random variable $W_i$.  At the first stage, it takes $n_{0}$ samples from alternative $i$ and computes the sample means and variances $(\bar x_{i,0}, \hat\sigma_{i,0}^2)$. By Bayes' rule, it indicates the prior distribution of $W_i\sim \mathrm{St}(\bar x_{i,0}, n_{0}/\hat\sigma_{i,0}^2,n_{0}-1)$, where $\mathrm{St}(\mu,\kappa,\nu)$ denotes the student-t distribution with mean $\mu$, precision $\kappa$, and degrees of freedom $\nu$. If additional $n_i-n_0$ samples are allocated to alternative $i$ and the overall sample mean and variance are $({\bar x}_i,{\hat\sigma}_i^2)$, then the posterior distribution $W_i$ becomes $\mathrm{St}({\bar x}_i, n_i/{\hat\sigma}_i^2,n_i-1)$. The final selection will go to the alternative with the largest sample mean, i.e., $(k)=\arg\max_{i} \bar x_{i}$. {A false selection will incur a linear loss which is $\max_i W_i -W_{(k)}$. Therefore, the problem for the second stage is to choose $(n_1,\cdots, n_k)$ to minimize the expected linear loss\footnote{In the original paper of \cite{chick2001new}, they also consider the sampling cost, which results an additional linear term in the objective function.}, i.e.,
\begin{align}\label{eqn: evi}
\min_{n_{1}+n_2+\dots+n_k=N} \mathrm E\left[\mathrm E\left[\max_i W_i -W_{(k)}\,\big|\,(X_{1}, X_{2}, \dots, X_{k})\right]\right],
\end{align}
Notice that $(n_1,n_2,\dots,n_k)$ are determined before the second stage. Therefore, to calculate the expected linear loss, we need to take the expectation with respect to the second-stage samples $(X_{1}, X_{2}, \dots, X_{k})$ which are random.} As the problem in \eqref{eqn: evi} has no closed-form solution, \cite{chick2001new} derived their allocation policy by asymptotically minimizing a bound of the expected loss.

{\cite{chick2001new} then adapted this two-stage procedure to the dynamic setting. At each intermediate stage $t$, a set of observations is collected from each alternative in the previous stages, based on which the current-stage allocation policy needs to be determined. This issue is essentially what the two-stage procedure above attempts to address. Therefore, it seems natural to determine this allocation policy by applying the two-stage procedure. Particularly, a myopic perspective is taken as if the selection is made at the end of stage $t$ and the current-stage allocation policy is then obtained by solving \eqref{eqn: evi},
\begin{align}
V_t^{\text{EVI}}(\mathcal E_t)=\min_{\sum_i n_{i,t}=\sum_i n_{i,t-1} +\tau} \mathrm E\left[\mathrm E\left[\max_i W_{i} -W_{{(k)}}\,\big|\,(X_{1,t}, X_{2,t}, \dots, X_{k,t})\right]\right],
\end{align}
where $X_{i,t}$ denotes the random samples that will be taken at stage $t$.} The extension from the above two-stage procedure to a sequential procedure encounters an obstacle, which is caused by the unbalanced samples from different alternatives. Technically, it involves the subtraction of two student-t random variables with different degrees of freedom. \cite{chick2001new} overcame this difficulty by using the Welch approximation.

\begin{algorithm}[ht]
\caption{EVI Procedure for Linear Loss}\label{procedure:EVI}
\vspace{3pt}
\textbf{Require:} Number of alternatives $k$, common first-stage sample size $n_0\geq 2$, total sampling budget $N$, sampling budget $\tau$ per stage.
\vspace{3pt}
\begin{algorithmic}[1]
    \State Generate $n_0$ samples from each alternative $i$.
    \State Set $t\leftarrow 0$, $n_{i,t}\leftarrow n_0$, $b_t\leftarrow\sum_{i=1}^k n_{i,t}$.
    \While{$b_t<N$}
    \vspace{3pt}
        \State Update $\bar{x}_i$ and $\hat\sigma_i^2$. Set $\bar x_{(1)}\leq \bar x_{(2)}\leq \dots\leq \bar x_{(k)}$ and $\mathcal{L}=\{1,2,\dots,k\}$.
        \vspace{3pt}
        \State Let $\lambda_{(i)(j)}^{-1} \leftarrow \hat\sigma_{(i)}^2/n_{(i),t}+\hat\sigma_{(j)}^2/n_{(j),t}$, $d_{(i)(k)}\leftarrow\bar{x}_{(k)}-\bar{x}_{(i)}$. Set $b_{t+1}\leftarrow b_t+\tau$.
        \vspace{3pt}
        \State For each alternative $(i)\in\mathcal L$, calculate
        \begin{small}
        \begin{align*}
		n_{(i),t+1} = \frac{\left(\tau+\sum_{(j)\in\mathcal L} n_{(j),t}\right)\left(\hat\sigma_{(i)}^2\eta_{(i)}\right)^{1/2}}{\sum_{(j)\in\mathcal L} \left(\hat\sigma_{(j)}^2\eta_{(j)}\right)^{1/2}}
        \end{align*}
        where
        \begin{equation*}
        \eta_{(i)} = \lambda_{(i)(k)}^{1/2}\frac{n_{i,t}-1+\lambda_{(i)(k)}d_{(i)(k)}^2}{n_{i,t}-2}\psi_{n_{i,t}-1}\left[\lambda_{(i)(k)}^{1/2}d_{(i)(k)}\right] \mbox{ for } (i)\neq (k), \mbox{ and } \eta_{(k)}=\sum_{(j)\neq (k)}\eta_{(j)},
        \end{equation*}
        and $\psi_s(\cdot)$ denotes the probability density function of a standard student-t distribution with $s$ degrees of freedom.
        \end{small}
        \While{$\min_{(i)\in \mathcal L}\,(n_{(i),t+1}-n_{(i),t})<0$ }
            \If{$n_{(i),t+1}-n_{(i),t}<0$}
            \vspace*{3pt}
                \State Set $\mathcal{L}\leftarrow \mathcal L\setminus(i)$ and $n_{(i),t+1} \leftarrow n_{(i),t}$.
            \EndIf
            \State For each alternative $(i)\in\mathcal L$, update
            \begin{small}
            \begin{equation*}
            \lambda_{(i)(k)}^{-1} = \left\{
            \begin{aligned}
            &\hat\sigma_{(i)}^2/n_{(i),t}+\hat\sigma_{(k)}^2/n_{(k),t}, &\mbox{if } (k)\in\mathcal L \\
            % &\hat\sigma_{(k)}^2/n_{(k),t}&\mbox{if } (i)\notin \mathcal L,(k)\in\mathcal L \\
            &\hat\sigma_{(i)}^2/n_{(i),t}, &\mbox{if } (k)\notin\mathcal L.
            \end{aligned}
            \right.
            \end{equation*}
            \end{small}
            \State Go back to Step 6.
        \EndWhile
        \State Generate $n_{i,t+1}-n_{i,t}$ samples from each alternative $i\in\mathcal L$. Set $t\leftarrow t+1$.
    \EndWhile
    \State Select $\arg\max \bar x_{i}$ as the best.
\end{algorithmic}
\end{algorithm}
The procedure is documented in Procedure~\ref{procedure:EVI}, where we assume the sampling cost from each alternative is the same and set as one. Notice that optimal sampling allocation policy in Step 4 looks very similar to that of the OCBA procedure (Step 5 in Procedure 3), because these two procedures are derived in a similar way as mentioned before.

In the same paper, \cite{chick2001new} also considered the problem with unconstrained budget, and proposed an EVI procedure to determine the number of replications to balance the replication costs against the reduction in expected opportunity cost. They also proposed analogous procedures with the 0-1 loss function. \cite{chick2010sequential} developed a variation of the EVI procedure in which the sampling budget is allocated to only one alternative at each stage. {For this special case, they showed that most of the approximations in solving this optimal allocation policy can be avoided and therefore derived a procedure with better performance, especially in the small-budget problems.}

%For the case when the sampling cost from each alternative is the same and unit one, we summarize the corresponding EVI procedure in Procedure 4. Notice that optimal sampling allocation policy in Step 4 looks very similar to that of the OCBA procedure (Step 5 in Procedure 3). This arises because these two procedures are derived in a similar way as mentioned before.

\subsection{One-step-look-ahead procedures}
%\begin{itemize}
%\item Multi-stage procedures and samples are taken one by one.
%\item One step look ahead and derive the myopic procedures.
%\item KG procedures
%\end{itemize}
%In this section we review the group of DP procedures which are featured with allocating samples one-by-one, which we name as the one-step-look-ahead procedures.
%Indeed, both OCBA and EVI have variations as one-step-look-ahead procedures.
In this section we review the group of DP procedures which are derived using the one-step look-ahead approximation. Specifically, we consider the knowledge-gradient (KG) procedure proposed by \cite{frazier2008knowledge}.

The KG procedure also adopts a Bayesian approach to solve the R\&S problem. Unlike the EVI, it determines the optimal sampling allocation policy by maximizing the expected terminal reward. Suppose that there is a budget of $N$ samples for selecting the best from $k$ alternatives.
Information collected from the samples is summarized in the posterior distribution of the unknown mean for each alternative. Let $\mu_i^t$ and $(\sigma_i^t)^2$ be the mean and variance of the posterior distribution for alternative $i$ after observing the first $t$ samples. Then, the problem is to determine the allocation of the $(t+1)$-th sample $z^{t+1}\in\{1,\cdots, k\}$ for $t=0,\cdots, N-1$, in order to
 maximize the expected terminal reward $\mathrm E\{\max_i \mu_i^{N}|\mathcal{E}_{N}\}$,  and the alternative with the largest $\mu_i^N$ is selected as the best. Here the information set $\mathcal E_t$ records the posterior mean and variance after the $t$th sample, and is updated according to the Bayes rule.

From the view of dynamic-programming formulation, we solve
%we can write the Bellman equation as
%\[
%V_t(\mathcal E_t) = \max_{z^{t+1}}\mathrm E\{V_{t+1}(\mathcal E_{t+1})|\mathcal E_t\},
%\]
%where
\begin{equation*}
V_t^\text{KG}(\mathcal E_t) = \max_{(z^{t+1},\cdots,z^{N})}\mathrm E\{V_{N}(\mathcal E_{N})|\mathcal E_t\}= \max_{(z^{t+1},\cdots,z^{N})}\mathrm E\{\max_i \mu_i^{N}|\mathcal E_t\}.
\end{equation*}
The key idea of KG procedure is to approximate $V_t(\mathcal E_t) $ by
\begin{equation*}
V_t^\text{KG}(\mathcal E_t) \approx \sum_{j=t}^{N-1} \max_{z^{j+1}}\mathrm E\{\max_i \mu_i^{j+1}-\max_i \mu_i^{j}|\mathcal E_j\} +\max_i \mu_i^t,
\end{equation*}
and the problem reduces to solve the one-step optimization problem
\begin{equation}\label{kg}
\max_{z^{t+1}}\mathrm E\{\max_i \mu_i^{t+1}-\max_i \mu_i^{t}|\mathcal E_t\}.
\end{equation}
Intuitively, it maximizes the increment (e.g., gradient) in the ``knowledge'' gained from the next sample, which explains the name ``knowledge gradient".

We assume that the samples across different alternatives are independent and have a common and known variance. In this special structure, the optimal solution of \eqref{kg} has a closed form (see Steps 4-5 in Procedure 5 or Theorem~1 in
\cite{frazier2008knowledge}). This structure is highly attractive from the implementational point of view. Besides, {the procedure possesses} other favorable properties. For instance, it is consistent, i.e., the selected alternative converges to the true best as the total sampling budget $N$ grows to the infinity, and {the suboptimality of the KG policy is bounded for any finite budget $N$}.

\begin{algorithm}[h]
\caption{KG Procedure}
\vspace{3pt}
\textbf{Require:} Number of alternatives $k$, total sampling budget $N$, common and known variance $\sigma^2$, prior predictive mean $\mu_i$ and variance $\sigma_i^2$ for each alternative.
\vspace{3pt}
\begin{algorithmic}[1]
    \State Set $t\leftarrow 0$. Let $\mu_i^t \leftarrow\mu_i, \beta_i^t\leftarrow 1/\sigma_i^2$ and $\beta=1/\sigma^2$.
    \While{$t< N$}
    \vspace{3pt}
        \State Calculate the variance of the change in predictive mean by taking a sample from alternative $i$, $\tilde\sigma_i^2=(\beta_i^t)^{-1}-(\beta_i^t+\beta)^{-1}$.
        \vspace{3pt}
        \State Calculate
        \begin{small}
        \begin{align}
        \zeta_i = -\bigg|\frac{\mu_i^t-\max_{j\neq i}\mu_j^t}{\tilde\sigma_i}\bigg|.
        \end{align}
        \end{small}
        \State Choose
        \[z^{t+1}=\arg\max_{i=1,2,\dots,k}\,\tilde\sigma_i (\zeta_i\Phi(\zeta_i)+\phi(\zeta_i)),\]
         {where $\Phi(\cdot)$ and $\phi(\cdot)$ denote the cumulative distribution function and probability density function of the standard Gaussian distribution, respectively.}
        \State Take a sample $y_{z^{t+1}}^{t+1}$ from alternative $z^{t+1}$. Update
        \[
        \beta_{z^{t+1}}^{t+1}\leftarrow \beta_{z^{t+1}}^{t}+\beta, \quad \mu_{z^{t+1}}^{t+1}\leftarrow (\beta_{z^{t+1}}^{t}\mu_{z^{t+1}}^{t}+\beta y_{z^{t+1}}^{t+1})/\beta_{z^{t+1}}^{t+1}.
        \]
        \State Set $t\leftarrow t+1$.
    \EndWhile
    \State Select $\arg\max \mu_{i}^t$ as the best.
\end{algorithmic}
\end{algorithm}
%The third family of procedures is the knowledge-gradient (KG) procedures. Unlike the OCBA and EVI procedures, the (KG) procedures determine the optimal sampling allocation policy to maximize the expected terminal reward and assume that the samples are allocated one-by-one, i.e., $\tau=1$. \cite{frazier2008knowledge} set the terminal reward as
%$V_T^*(\mathcal \xi)=\max_i \mathrm E\{U_i|\mathcal{E}_{T}\}=\max_i \mu_i^T$ and solve the approximated Bellman equation \eqref{eqn:bellman_myopic} as
%\begin{align}\label{bellman_kg}
%\arg\max_{n_{1,t}+n_{2,t}+\dots+n_{k,t}=t}\ \mathrm E\{\max \mu_i^{t+1}|\mathcal{E}_{t}\}.
%\end{align}
%For the VFA, the KG procedure substitutes the value function by $(\max \mu_i^{t+1}-\max \mu_i^t)$, which refers to the increment in the value of information by taking an additional sample.  Interestingly, the procedure turns to maximize the increment in the value of information in order to obtain the optimal allocation policy, i.e.,
%%and the corresponding optimization ]\eqref{bellman_kg} is written as
%\begin{align}\label{bellman_kg1}
%\arg\max_{n_{1,t}+n_{2,t}+\dots+n_{k,t}=t}\ \mathrm E\{\max \mu_i^{t+1}-\max \mu_i^t|\mathcal{E}_{t}\}.
%\end{align}
%Notice that \eqref{bellman_kg} and \eqref{bellman_kg1} are equivalent since $\max \mu_i^t$ is known under $\mathcal{E}_t$.

The original KG procedure of \cite{frazier2008knowledge} has several variations. For instance, \cite{frazier2009knowledge} and \cite{Xie2016} extended the procedure to the case of correlated sampling and correlated Gaussian beliefs on the mean vectors. \cite{Ryzhov2016} {adopted} a different way to define the value of information functions and then derived the corresponding optimal sampling allocation rule. In this rule, the allocation ratios among the non-best alternatives are quite similar to that of the OCBA procedures. However, the total proportion of samples allocated to these non-best alternatives vanishes as the total sampling budget grows to infinity.
To understand the connection between Ryzhov's procedure and the OCBA procedure, \cite{PengFu2017} showed that the allocation rules of the OCBA procedure can be achieved by slightly modifying the function used to describe the value of information in \cite{Ryzhov2016}.

\section{Large-scale R\&S procedures using parallel computing}

As mentioned before, many existing R\&S procedures, under either the fixed-precision or the fixed-budget formulations, are designed to solve small- or medium-scale problems, with total number of alternatives typically less than $500$, which is largely  due to the limited computing resource.  On one hand, there are many large-scale R\&S problems in practice that have thousands to millions of alternatives, which are traditionally solved by  optimization-via-simulation (OvS) algorithms (see, for instance, \cite{HongNelson2009} and \cite{HNX2015} for comprehensive reviews of OvS). On the other hand, the fast development of computer technology and parallel computing (e.g., either from the multi-core personal computers to many-core servers or from smart phones to cloud services) are prevalent and ready for ordinary users to access. Then, using parallel computing to solve a large-scale R\&S problem directly becomes an interesting research topic. It has even been labeled as one of the three central developments in the past 15 years by \cite{FuHenderson2017}.

Researchers begin investigating parallel computing for R\&S problems by asking the following questions: (i) Can existing R\&S procedures can be easily implemented in a parallel fashion; (ii) If not, how are these procedures modified to suit for parallel computing environments? (iii) In the process of parallelization, what kind of substantial issues need to be addressed?
To the best of our knowledge, \cite{YLCL2001} and \cite{Chen2005} are the two earliest works in the literature that try to answer the first question.
In particular, the former implemented the OCBA procedure in a web-based parallel environment, and the latter executed a multi-stage procedure by
distributing the simulation tasks to multiple processors. However, both studies tested their procedures only for a small-scale problem with 10 alternatives, so it is not clear whether their procedures are suitable for handling large-scale problems.  \cite{LHNW2015} (and their conference paper
\cite{LuoHong2011}) and \cite{NCHH2017} (and their conference papers  \cite{NHH2013,NHH2014,NCHH2015}) are works that intend to answer the three questions. They demonstrated that redesigned procedures can be used to solve large-scale problems with thousands to millions of alternatives in different parallel computing environments.

There are various parallel computing environments that are suitable for R\&S problems, and they can be in general classified into three categories, i.e., Message-Passing Interface (MPI), Hadoop MapReduce and Apache Spark (\cite{NCHH2017}). The MPI \citep{gropp1999using} is a standardized and portable message-passing protocol for parallel programming on several parallel computing architectures, which is equipped with C/C++ and Fortran libraries.
Hadoop \citep{dean2008mapreduce} is an open-source framework designed for distributed storage and processing of large amounts of data and computation using the MapReduce programming architecture.  Apache Spark \citep{zaharia2010spark}  is also an open-source framework for general-purpose parallel computing.
Both MapReduce and Spark are supported by several commercial clouds including Amazon EC2,  Google Cloud Platform and Microsoft Azure (\cite{ZhongHong2019}). Note that all the three frameworks can be implemented using the Master/Worker parallel structure. MPI allows more flexibility of parallel implementation but does not detect or manage core failures automatically compared with MapReduce and Spark.

In the following, we first briefly describe both the theoretical and implementational challenges  as modifying existing R\&S procedures to suit for parallel computing environments in Section~\ref{sec:parallelproc}. Then, we introduce some different performance measures and new frameworks that are developed for large-scale R\&S problems in Section~\ref{sec:knockout}. Two representative procedures are also presented in Sections~\ref{sec:parallelproc} and \ref{sec:knockout}, respectively.

\subsection{Extending existing procedures to parallel}
\label{sec:parallelproc}

In traditional R\&S problems, the efficiency of a procedure can often be measured by the total running time, which is approximately the total simulation time of generating observations from different alternatives. This method is reasonable since the operations of all other calculations and comparisons  are quite fast and the total time of these operations could be negligible compared with the total simulation time as solving small-scale problems  in a single-processor environment. However, when handling large-scale problems in a parallel computing environment, the situation becomes complicated since the comparison operations may become the bottleneck. The  communications and synchronizations among different processors may also need to be taken into consideration. In other words, to measure the efficiency of a procedure in parallel computing environments, we shall evaluate the running time from four aspects, i.e., the simulation time, the comparison time, the communication time  (i.e., the time to transfer information between different processors) and the synchronization time (i.e., the time to wait for the ready state of all processors).
For the sake of presentation,
we take the stage-wise and fully sequential procedures in the fixed-precision formulation to illustrate the tradeoff among the four aspects.

Stage-wise procedures are easy to parallelize and there are no communications among processors until the comparison operation at the end of each stage, which means they are efficient in comparison and  communication. The synchronization is also not an issue if the simulation tasks are distributed evenly onto different processors.
Compared with fully sequential procedures, however, stage-wise procedures are typically not efficient in total sample size, i.e., inefficient in simulation.
For fully sequential procedures,   they conduct all-pairwise comparisons (i.e., $k(k-1)/2$ in the worst-case) among all alternatives still in contention at each round when all alternatives add one observation, implying frequent communications and synchronizations among different processors. Therefore, they are inefficient in comparison, communication and synchronization.

Note that the total sample size is inherently determined by the theoretical framework of a procedure, which could be hardly  reduced even using parallel computing. Therefore, there is little room for improving the efficiency of stage-wise procedures, and  many works in the literature focus on improving the efficiency of fully sequential procedures by redesigning them to be fit for parallel computing in order to
reduce the times for comparison, communication and synchronization. For instance, \cite{LHNW2015} address the synchronization issue by
proposing an asynchronization scheme to achieve a high simulation efficiency of sampling, and point out the potential issues caused by all-pairwise comparisons and frequent communications. Later,  \cite{NCHH2017} and \cite{ZLLH2019} addressed the comparison issue using two different approaches, namely, a ``divide-and-conquer'' scheme by distributing the all-pairwise comparisons and a new comparison scheme by defining the ``best'' alternative differently. They further mitigated the communication burden by using batching techniques and boosting the sample size of all surviving alternatives to a maximum number afterwards. Notably, different batching techniques and boosting methods are proposed in \cite{NCHH2017} and \cite{ZLLH2019}, thus resulting in different theoretical foundations of their procedures.
Before introducing more details, we first briefly describe the aforementioned Master/Worker parallel structure which has been used in \cite{LHNW2015}, \cite{NCHH2017} and \cite{ZLLH2019}.

Suppose that there are $m+1$ processors in the parallel computing environment, in which one processor serves as the master and the rest $m$ processors serve as the workers, denoted by workers $1,2,\ldots,m$. The master is the controller who determines the start and stop of the program, creates $m$ job tasks for the workers, manages the data information and performs all other necessary calculations. Workers  $1,2,\ldots,m$ function in a simple way: Taking the  task  from the master, processing the task, submitting the result to the master and requesting the next task. The communications occur only between the master and workers, and there is no communication among workers.

To address the synchronization issue,
\cite{LHNW2015} defined each job task as generating one observation from one alternative that is still in contention, and all alternatives in contention are queued in a round-robin order in front of the master. This one-by-one task assignment scheme requires no synchronization among workers, and can indeed fairly balance the workloads of different workers.
However, given the random processing time of each task on different workers, the sequence of the simulation results sent back to the master is also random, which is likely to be different from the round-robin assigning order. Therefore, it may cause unexpected statistical issues as implementing existing fully sequential procedures based on the output sequence.
One straightforward way is to restore all the outputs exactly in the same order as in the original input sequence, and perform comparisons according to the input sequence, which is the basic idea for the vector-filling $\mathcal{KN}$ (VKN) procedure of \cite{LHNW2015}. We omit the details about the VKN procedure, but summarize two disadvantages of VKN. First, it needs to create a vector to store these outputs, which  may require a large amount of memory. Second, it does not allow the failure of any worker or the communication interruption, as the vector may be incomplete if some of the simulation results are lost in these situations.

To resolve these problems, \cite{LHNW2015} proposed the asymptotic parallel selection (APS) procedure, which performs all-pairwise comparisons directly based on the output sequence and introduces a phantom alternative to determine the time points for comparisons.\footnote{Note that the phantom alternative does not need to be processed by any worker but immediately returns to the {master} for requiring the master to perform  eliminations by  conducting all-pairwise comparisons.} However, the desired finite-time PCS guarantee is no longer achieved, as the sample sizes of different alternatives in the output sequence are random and perhaps unequal; these observations from the same alternatives are not even i.i.d. Fortunately, the innovative idea of introducing the phantom alternative, which serves as a  drumbeat process with predetermined time points, allows to establish a finite lag of the difference between the input and output sequences that finally vanishes in an asymptotic regime. By doing so, the APS procedure of \cite{LHNW2015} provides an asymptotic PCS guarantee. We present  the APS procedure in Procedure~\ref{proc:APS}.

\begin{algorithm}[htp]
\caption{Asymptotic Parallel Selection (APS) Procedure}\label{proc:APS}
\vspace{3pt}
\textbf{Require:} Number of alternatives $k$, PCS $1-\alpha \in (1/k,1)$, IZ parameter $\delta >0$, the first-stage sample size $n_0\geq 2$, and the number of processors $m+1$  (i.e., one master and $m$ workers).
\vspace{3pt}
\begin{algorithmic}[1]
    \State Let $a=-\log \left[2\alpha/(k-1)\right]$ and $I\leftarrow \{1,2,\dots,k\}$.
    \State Let $p$ denote the phantom alternative queued after each round-robin cycle and let the stage $r$ denote the current sample size of $p$.
    \State For all $i\in I$, record the triple $\big(N_{ir},\sum_{\ell=1}^{N_{ir}}Y_{i\ell},\sum_{\ell=1}^{N_{ir}}Y_{i\ell}^2\big)$, where  $Y_{i\ell}$ is
the $\ell$th completed observation from alternative $i$ and $N_{ir}$ is the total sample size obtained from alternative $i$ at stage $r$.  Set $r \leftarrow 1$ and set $N_{ir} \leftarrow 0$, $\sum_{\ell=1}^{N_{ir}}Y_{i\ell} \leftarrow 0$, $\sum_{\ell=1}^{N_{ir}}Y_{i\ell}^2\leftarrow 0$ for all $i\in I$.
    \State Using the round-robin order to start simulations on workers $1,2,\ldots,m$.

    \While{$|I|>1$}
    \State Wait for the next observation $Y_{h\cdot}$.
	\If {$h\in I$}
      	 \State Update  $\sum_{\ell=1}^{N_{hr}}Y_{h\ell} \leftarrow \sum_{\ell=1}^{N_{hr}}Y_{h\ell} + Y_{h\cdot}$, $\sum_{\ell=1}^{N_{hr}}Y_{h\ell}^2\leftarrow \sum_{\ell=1}^{N_{hr}}Y_{h\ell}^2 +Y_{h\cdot}^2$ and $N_{hr} \leftarrow N_{hr}+1$.%, and break;
      	 	\ElsIf{$h\notin I$  and $h \neq p$}
    	 \State Drop the observation.%, and break;
	\ElsIf{$h = p$}
%    	  \State Calculate the sample mean and sample variance for all $i\in I$
%    \begin{footnotesize}
%    \begin{eqnarray*}
%\bar{Y}_i(N_{ir})= \frac{1}{N_{ir}}\sum_{\ell=1}^{N_{ir}}Y_{i\ell}, \mbox{ and }
%S_i^2(N_{ir})    = \frac{1}{N_{ir}-1} \left[\sum_{\ell=1}^{N_{ir}}Y_{i\ell}^2
%       -\frac{1}{N_{ir}}\left(\sum_{\ell=1}^{N_{ir}}Y_{i\ell}\right)^2\right].
%\end{eqnarray*}
%\end{footnotesize}
  \ForAll{ $i,j\in I$ and $i\neq j$}
\If {$N_{ir}\geq n_0$ and $N_{jr}\geq n_0$}
    \State Let
   \begin{footnotesize}
\begin{eqnarray*}
\tau_{ij,r}=\left[\frac{S_i^2(N_{ir})}{N_{ir}}+\frac{S_j^2(N_{jr})}{N_{jr}}\right]^{-1},
\end{eqnarray*}
\end{footnotesize}
where $S_i^2(N_{ir})$ and $S_j^2(N_{jr})$ are the sample variance of alternatives $i$ and $j$, respectively.
\Else
        \State  Let $\tau_{ij,r}=0$.
\EndIf
\EndFor

 \State Set
 \begin{small}
\[
~~~~~~I \leftarrow \left\{i: i\in I \mbox{ and }
\tau_{ij,r}\left[\bar{Y}_i(N_{ir}) - \bar{Y}_j(N_{jr})\right]
\geq \min\left\{0,-\frac{a}{\delta}+\frac{\delta}{2}\tau_{ij,r}\right\},
\forall j \in I, j\neq i \right\},
\]
 \end{small}
 where $\bar{Y}_i(N_{ir})$ and $\bar{Y}_j(N_{jr})$ are the sample mean of alternatives $i$ and $j$, respectively.

\State Set $r\leftarrow r+1$.
	 \EndIf

    \EndWhile
    \State Select the alternative in $I$ as the best.
\end{algorithmic}
\end{algorithm}

As mentioned but not addressed by \cite{LHNW2015},  all-pairwise comparisons conducted in the master could overwhelm the workload of the master and the frequent communications between the master and workers could become the bottleneck for solving large-scale R\&S problems.
In order to reduce the comparisons,
the good selection procedure (GSP) of \cite{NCHH2017} proposes a ``divide-and-conquer'' approach to distributing  all-pairwise comparisons onto the workers. The master initially divides $k$ alternatives into $m$ groups and asks each worker to conduct local all-pairwise comparisons for eliminations within the assigned group. Then, the computational complexity of comparisons at each round is reduced from $\mathcal{O}(k^2)$ to $\mathcal{O}(k^2/m^2)$.
To further improve the elimination efficiency,  at the beginning of each local comparison round, the master retrieves the $m$ local bests  from the $m$ groups to find the global best and then sends the global best to the $m$ groups for additional comparisons.

To reduce the
frequent communications, GSP introduces a batching technique of samples. In particular, it suggests each worker  to simulate a batch of 100 or 200 samples from one alternative once at a time. In addition, when a surviving alternative {takes} enough samples, i.e., reaching a threshold, the procedure requires all surviving alternatives take a maximum number of samples. It {selects} the one with the largest sample mean as the best. By doing so, GSP can significantly reduce the communication frequency. In fact, the maximum sample size in the boosting stage (i.e., Stage 3 in their procedure) is constructed based on the Rinott's result and the sequential elimination rule in Stage 2 is built on the results of
\cite{Hong2006}. Taking advantage of both sequential and stage-wise frameworks, GSP is an excellent hybrid procedure that  not only  improves the comparison and communication efficiency, but also provides a finite-time PGS guarantee of \eqref{PGS}.

One potential drawback of GSP is its conservativeness in terms of total sample size due to the batching technique of samples and the error separation in the hybrid structure. In other words, GSP sacrifices a certain level of sampling efficiency to achieve lower computational complexities of comparisons and communications as well as the finite-time statistical guarantee.  This drawback motivated \cite{ZLLH2019} to design the parallel Paulson's procedure (PPP) that takes a different approach to achieving simulation, comparison and communication efficiency.

In terms of the simulation efficiency, PPP  adopts
 the well-known sequential procedure of \cite{paulson1964sequential}, which often requires a smaller sample size than the stage-wise and the hybrid procedures in the same desired guarantee level. In terms of the comparison efficiency, PPP breaks all-pairwise comparisons into comparisons with the ``best'', which  reduces the computational complexity of comparisons from $\mathcal{O}(k^2)$ to $\mathcal{O}(k)$ at each comparison round. The ``best'' involves both the sample mean and sample  variance information, which is different from the global best involving only the mean information defined in GSP.
 In terms of the communication efficiency, PPP uses a different batching technique, i.e., batching alternatives instead of  samples, that can reduce the communication frequency without increasing the sample size. PPP also allows to boost the sample size to a maximum number in Paulson's sequential framework. In addition, by incorporating the result
of \cite{kao1980sequential}, PPP can also achieve the PGS guarantee.

We omit the presentation of both GSP and PPP. Interested readers may refer to the works of the abovementioned authors for more details. The CRNs technique is difficult to implement for VKN, APS, and GSP because of the asynchronized simulation scheme in \cite{LHNW2015} and \cite{NCHH2017}. Moreover, CRNs are not suitable for PPP because of the newly designed comparison scheme in \cite{ZLLH2019}.

\subsection{New parallel framework with sample-size optimality}
\label{sec:knockout}	

The aforementioned procedures for parallel computing environments are all built on the paradigm of existing stage-wise or fully sequential R\&S procedures. They have even successfully solved R\&S problems with up to millions of alternatives. Therefore, we must ask whether they are fundamentally suitable for handling large-scale problems. More precisely, we would like to know how the expected total sample size $\mathrm  E[N]$ increases as the number of alternatives $k$ increases.

To answer the question, \cite{ZhongHong2019} first proved that {the growth rate} of $\mathrm  E[N]$ for any procedure with the PCS guarantee is lower bounded by  $\mathcal{O} (k)$, and then defined the sample-size optimality of a procedure if the upper bound of the growth rate of $E[N]$ can achieve the lower bound rate, that is  $\mathrm E [N] = \mathcal{O} (k)$. Intuitively speaking, the lower bound of $\mathrm E[N]$ is easy to understand since each alternative requires at least one observation to estimate the unknown  mean, resulting in a total sample size growing at least linearly in $k$. The lower bound is universal for all stage-wise and fully sequential procedures in either the IZ or IZ-free framework \citep{ZhongHong2019}.

However, the upper bound is typically higher than the order of $k$ for all existing stage-wise and fully sequential procedures. For instance, the expected total sample size of each alternative in Paulson's procedure grows
proportionally to the ending point of the continuation region, which grows in the order of $\log k$, leading to $\mathcal{O}(k\log k)$ in total sample size. This is because it needs to compare the best with all other $k-1$ alternatives in pairs in theoretical formulation,
 which is also true for other fully sequential procedures such as $\mathcal{KN}$ procedure. Similarly, in stage-wise procedures, e.g., one-stage procedure of \cite{bechhofer1954single} and two-stage procedure of \cite{rinott1978two}, they also need to compare the best with all other $k-1$ alternatives in a joint formulation as in \eqref{Bechhofer_typeI}, so the sample size of each alternative in  \eqref{bechhofer_n} grows as $k$ increases, which inevitably leads to a higher order of $k$ in total sample size. In other words, neither fully sequential nor stage-wise procedures can achieve the sample-size optimality due to the nature of comparisons between the best alternative and all others.

Inspired by the  knockout tournament arrangement of tennis Grand Slam tournaments, \cite{ZhongHong2019} proposed a novel parallel selection framework in which the champion (i.e., the unknown best) does not have to play with all others to be declared as the best. In particular, the knockout tournament ($\mathcal{KT}$) procedures of \cite{ZhongHong2019} divide all alternatives into pairs and construct a ``match'' between two alternatives in pair, keeping the winner for the next round ``matches'' while knocking out the loser after the current round. By doing so,  $\mathcal{KT}$ procedures  eliminate approximately half of the alternatives at each round and therefore achieve the theoretical lower bound of $\mathrm E[N]$. The sample-size optimality is achieved regardless of whether the variances of the alternatives are known or not. However, it will only change the constants on the optimal upper bounds.

This structure of $\mathcal{KT}$ procedures is perfect for parallel computing in terms of synchronization,  communication and comparison efficiency.
$\mathcal{KT}$ procedures can divide all alternatives into $m$ subsets, and assign each subset to one processor to conduct selections simultaneously among the alternatives in that subset. Then, neither synchronization nor communication  are necessary among different processors until each processor produces a local best alternative.
Since all ``matches'' in each processor are conducted independently and locally, and CRNs technique can be easily implemented into  $\mathcal{KT}$ procedures. In addition, because comparisons are only made within ``matches'', \cite{ZhongHong2019} demonstrated that  the comparison time in the procedures is negligible compared with the simulation time. In fact, the number of comparisons is only half of the total sample size since simulating two observations (i.e., one for each alternative) is coupled with just a single comparison.

In each ``match'', any existing R\&S procedures can be used to determine the winner, and $\mathcal{KT}$ procedures adopt $\mathcal{KN}$ procedure to achieve  the PCS guarantee and to gain sampling efficiency on total sample size. The sampling efficiency can be further improved by assigning more than two alternatives into one ``match'',  and the PGS guarantee can also be obtained by allocating the IZ parameter in different rounds.

\begin{algorithm}[htp]
\caption{ {\boldmath$\mathcal{KT}^+$} Procedure }\label{proc:KTP}
\vspace{3pt}
\textbf{Require:} Number of alternatives $k$, PCS $1-\alpha$ $\left(0<\alpha<1-1/k\right)$, IZ parameter $\delta >0$, the first-stage sample size $n_0 >2$, the  parameter $\lambda$ ($0<\lambda<\delta$),  number of alternatives $g\ge2$ within a ``match'', and the number of processors $m+1$  (i.e., one master and $m$ workers).
\vspace{3pt}
\begin{algorithmic}[1]
		\State Let $I_r^s$ be the set of alternatives in contention at the beginning of round $r$ in processor $s$ for $s=1,2,\ldots,m$.
		\State Equally allocate $k$ alternatives to $m$ processors so that each processor handles the selection of approximately $ k/m$ alternatives, e.g., for $i=1,2,\ldots,k$, let,
		\begin{equation*}
		I_1^{(i\bmod m)+1} = I_1^{(i\bmod m)+1}\cup\left\lbrace i\right\rbrace.
		\end{equation*}
		\ForAll{$s=1,2,\ldots,m$}
		\State Set $r=1$.
		\While{$\left|I_r^s\right|>1$}
		\State Let $I_{r+1}^s=\emptyset$. Group alternatives in $I_r^s$ with the size of $g$. In case of leftover ones, let them form a group. After grouping, there are in total $\left\lceil\left|I_r^s\right|/g\right\rceil$ groups. Let $I_{r,q}^s$ denote the set of the alternatives in group $q$ for $q=1,2,\ldots,\left\lceil\left|I_r^s\right|/g\right\rceil$ of processor $s$ at round $r$.
		\State Let $\alpha_r=\alpha/2^r$. For each group $q=1,2,\ldots,\left\lceil\left|I^r_s\right|/g\right\rceil$, set $\mathcal{C}=I_{r,q}^s$ and compute,
			\begin{equation*}\label{eqn:KTMATCH}
			I_{r+1}^s=I_{r+1}^s\cup\big\lbrace \mathcal{KN}\left(\mathcal{C},\alpha_r,\delta,n_0\right)\big\rbrace.
			\end{equation*}	
		\State Set $r=r+1$.
		\EndWhile
		\State  Let $I_s$ denote the index of the alternative in $I_r^s$.
		\State Take $n_0$ observations from alternative $I_s$. Calculate its sample variance $S^2_{I_s}$ based on these $n_0$ observations. Set $r=\left\lceil \log_g\frac{k}{m}\right\rceil+1$, $\alpha_r=\alpha/2^r$, and $h\left(\alpha_r,m,n_0\right)$, where $h\left(\alpha_r,m,n_0\right)$ is the Rinott's constant determined by $\alpha_r$, $m$, and $n_0$.
		\State Set
		\begin{equation*}\label{eqn:KTMATCH2}
			N_{\max, I_s}= \max\left\lbrace n_0,\left\lceil\left(\frac{h\left(\alpha_r,m,n_0\right)S_{I_s}}{\delta}\right)^2\right\rceil\right\rbrace,
			\end{equation*}
			Then, take additional $N_{\max, I_s}-n_0$
			observations from alternative $I_s$.
		\State Compute the sample mean $\bar{X}_{I_s}(N_{\max, I_s})$.
		\EndFor
		\State Let $I$ denote the set of alternatives containing all the best alternatives produced by $m$ processors. Select the alternative with the largest sample mean $\bar{X}_{I_s}(N_{\max, I_s})$ for $I_s \in I$.
	\end{algorithmic}
	\end{algorithm}

For the simplicity of presentation, we adopt the notation of $\mathcal{KN}\left(\mathcal{C},\alpha_r,\delta,n_0\right)$ in \cite{ZhongHong2019} to denote the output of executing the $\mathcal{KN}$ procedure in each ``match'', which provides a PCS of  $1-\alpha_r$ among the alternatives with unknown means and unknown variances in the set $\mathcal{C}$ when the first stage sample size is $n_0$ and the IZ parameter is $\delta$. In the following, we present the $\mathcal{KT}^+$, one of the $\mathcal{KT}$ procedures for parallel computing environments of \cite{ZhongHong2019}, in Procedure~\ref{proc:KTP}.

We conclude this section by briefly reviewing some other recent works on parallel R\&S problems.
\cite{HunterNelson2017} argued that  different performance measures and different formulations are needed for large-scale R\&S problems. As a response,
\cite{pei2018new} proposed a different objective function, i.e., the expected false elimination rate (EFER) rather than the PCS/PGS. They argued that having a sizeable number of alternatives is more reasonable. \cite{ma2018sequential}  extended the Envelop procedure of \cite{ma2017an}
to parallel computing environments, which established a new error bound of Brownian motion processes inspired by the multi-armed bandit problem of
 \cite{Jamieson2014}.
Notice that these abovementioned procedures using parallel computing are fixed-precision procedures. Some of the well-known fixed-budget R\&S procedures  have also been adapted to parallel computing environments. For instance,
 \cite{KP2018} considered the asynchronization issue as extending both OCBA and KG procedures in parallel computing environments and discussed the efficiency of these procedures for small- and large-scale problems.

\section{Emerging R\&S problems}

Besides fixed-precision and fixed-budget R\&S problems, some research problems that expend classical R\&S from different perspectives have emerged, e.g., considering multiple performance measures, taking the input uncertainty into account, and treating the performance measure as a function of the underlying contexts. In the following, we  briefly discuss recent achievements in these topics, without presenting the detailed procedures.

\subsection{Constrained R\&S}

Traditional R\&S problems often focus on only one
performance measure. However, in various practical situations, we may be interested in multiple performance measures. For instance, in service centers managers  are concerned about the expected cost and customer waiting times. In
production systems managers care about not only the expected throughput  but also the associated product quality.
One way to deal with multiple performance measures is to model the primary one as the objective and to model others as constraints. This leads to
constrained R\&S problems considered in the simulation literature.

The study by \cite{AndradottirKim2010} is one of the first works in this area. The authors indicated that the primary and secondary simulation outputs are i.i.d. bi-variate Gaussian random vectors with unknown mean vector and covariance matrix. They extended the IZ formulation to this problem and designed fixed-precision procedures that are capable of solving the problem. \cite{HAK2014} further developed the work of \cite{AndradottirKim2010} to handle multiple secondary performance measures. Then, \cite{healey2015minimal} reconsidered the problem by taking the switching cost into consideration. Rather than modeling the secondary performance as a Gaussian random vector as all previous works, \cite{HLN2015} perceived the secondary performance as a Bernoulli random variable and the secondary performance measure as a probability. Hence, they called the problem a chance-constrained R\&S problem. They built a hypothesis test on the chance constraint, thus resulting in an efficient two-stage procedure that performs the feasibility check in the first stage and selects the best among all the sample feasible alternatives in the second stage.

Different from the IZ formulation of the constrained R\&S mentioned above, \cite{lee2012approximate} addressed the problem under the OCBA framework. Meanwhile, \cite{HunterPasupathy2013} and \cite{pasupathy2015stochastically} solved the problem based on the large deviations theory by analyzing the asymptotic rate of identifying the optimal feasible solution, which is later known as Sampling Criteria for Optimization using Rate Estimators (SCORE) simulation framework. Recently, \cite{GGXS2018} considered constrained R\&S problems in the OCBA formulation, in which they used large deviations theory and incorporated a quadratic regression metamodel to improve the efficiency further. %Taking advantage of the IZ-free formulation, \cite{CFLT2020} demonstrated the benefit of designing a fully sequential procedure that simultaneously conducts feasibility checking and optimality checking (i.e., comparisons among the secondary and primary performance measures of different alternatives, respectively).

\subsection{Multi-objective R\&S}

Except for the constrained R\&S formulation, another way to handle multiple performance measures is to treat them as simultaneous objectives to optimize, thereby giving rise to the multi-objective R\&S formulation. Multi-objective R\&S problems {are different from} the classic single-objective R\&S problems mainly in two aspects. First, in the multi-objective problems, that
an alternative ``dominating'' another alternative means that it is better on all objectives. Therefore, the single ``best" alternative that dominates all  others may not exist, and  the goal of multi-objective R\&S turns to select the set of  all non-dominated alternatives, termed by the Pareto set. Second, two types of probabilities are defined to measure the errors when a dominated alternative is included into the Pareto set and a non-dominated alternative is excluded from the set.

Multi-objective R\&S procedures are often designed based on traditional R\&S procedures. Interested readers may see \cite{hunter2019an} for a detailed review. \cite{lee2010finding} extended the OCBA procedures to find the sampling allocation rule that minimizes the weighted sum of the two types of error probabilities. Alternatively, \cite{FeldmanHunter2018} and \cite{AHP2020} derived the optimal sampling allocation rule by adopting the SCORE framework. \cite{branke2015a} supplemented the EVI procedure. At each stage, the new procedure allocates samples to the alternative that most likely changes the observed Pareto set. \cite{branke2016multiobjective} further developed the KG procedures and allocated the sample at each stage such that the estimated expected Pareto set is the closest to the true Pareto set.

The multi-objective procedures above are all designed for fixed budget. Meanwhile, the literature has also suggested several fixed-precision multi-objective procedures. \cite{batur2018methods} formulated the mean-variance portfolio analysis as a bi-objective R\&S problem and propose a bi-objective procedure under the IZ formulation that controls both types of error probabilities. \cite{wang2017sequential} designed a sequential IZ-free procedure for the multi-objective procedure based on the generalized sequential likelihood ratio test.

\subsection{R\&S with input uncertainty}

In simulation studies, the input distributions are often estimated from the data and other information, thus having uncertainty in them. This uncertainty is called input uncertainty. For instance, when modeling the arrival process of an online service system, multiple plausible distributions can fit the input historical data, especially when the data set is not sufficiently large. When specifying the demand curve of a newly launched product, different managers may have varying beliefs on it. Any distribution of such items can be used as a proxy of the true input distribution. However, the corresponding “best” alternative may be different. In other words, no single alternative may be the best for all the possible scenarios of the input distributions. Then, taking the uncertainty of the simulation model into consideration in making R\&S decisions is an interesting problem in R\&S.

\cite{song2015input} studied the impact of input uncertainty on the classic IZ procedures and found that a straightforward application of IZ procedures may fail to deliver the desired PCS in the presence of uncertainty. They further proposed an adjustment to provide an average PCS. However, this average PCS guarantee cannot be delivered for some configurations of the competing alternatives. Therefore, it is still necessary to design procedures that are able to address the R\&S with input uncertainty.

To design such procedures, \cite{fan2013robust} innovatively proposed a robust selection-of-the-best (RSB) framework. Particularly, the RSB formulation includes all the possible scenarios of input distribution into an ambiguity set and then takes a robust perspective to define the best alternative with respect to the worst-case mean performance measures over the ambiguity set. \cite{fan2019distributionally} further improved their work and proposed both two-stage and sequential procedures that can achieve the user-specified PCS. These RSB procedures are also tested by a healthcare queueing system with both synthetic and real hospital data. \cite{gao2017robust} considered the RSB formulation under the OCBA framework, in which the approximated PCS is also measured by the worst-case performance. Besides the PCS guarantee, \cite{WuZhou2019} took the RSB formulation from the fixed-budget viewpoint into account, in which a joint budget for both collecting input data and running simulations is given in advance.

\subsection{R\&S with covariates}

In some problems, the performance measure of an alternative may vary as a function of some
observable random covariates, which are also known as side information, auxiliary quantities,
or contextual variables. For instance, in healthcare management, the treatment outcome of one particular drug may depend on the patient's biometric characteristics. In revenue management, the best assortment could vary according to customer segmentation. Then, the best alternative is not universal but depends on the value of underlying covariates (e.g., patient's biometric characteristics or customer segmentation), and this type of selection of the best problem is called R\&S with covariates (R\&S-C) or contextual R\&S.

Reasonably defining and measuring a correct selection of the best is the first question that needs to be addressed. \cite{shen2017ranking} were the first to introduces several definitions of correct selection for R\&S-C from the frequentist perspective. They first defined the {\it conditional} PCS, which is denoted by PCS($\bf x$), as the probability of selecting the best alternative (more precisely, the good alternative within IZ) for an individual whose random covariates (denoted by $\bf X$) take the value $\bf x$. Then two forms of {\it unconditional} PCS are introduced. One is the average PCS, i.e., $\mathrm E[{\rm PCS}({\bf X})]$. The other is the worst PCS, i.e., 
$\min_{x\in \Omega} {\rm PCS}({\bf x})$
, where $\Omega$ is the support of $\bf X$. In both \cite{shen2017ranking}  and the subsequent work  \cite{SHZ2019}, they assumed a linear model between the mean performance of an alternative and the corresponding covariates. They developed fixed-precision procedures that can produce selection policies (mapping from covariates to alternative index) to achieve the desired targets of unconditional PCS. The IZ formulation in R\&S-C is natural and critical, since the mean performance surfaces of alternatives may intersect somewhere and the performance of different alternatives at the neighborhood of intersection points can be arbitrarily close or equal. \cite{li2018data} adopted the R\&S-C framework developed by \cite{shen2017ranking}. However they designed new selection procedures to accommodate the high-dimensional covariates and the general (nonlinear) dependence between the mean performance of alternatives and the covariates.

Fixed-budget R\&S-C problems have also been tackled under the Bayesian framework, with the aim of adaptively allocating a given sampling budget to the alternatives and over the domain of covariates to find the best response across the entire domain of covariates efficiently. \cite{hu2017sequential} proposed to model the performance functions of alternatives as Gaussian random fields and used the expected improvement criteria to develop Bayesian procedures. \cite{pearce2017efficient} followed the same setting in Hu and \cite{hu2017sequential} and proposed a KG-based sampling policy with a focus on efficiently estimating the expected improvement over a continuous domain. \cite{ZSHD2020} extended the problem to a more general setting where the sampling noise can be heteroscedastic and the sampling cost at different locations can be different. More importantly, they provided a theoretical analysis of the asymptotic behavior of the KG based policy and proved that the best alternative as a function of the covariates will be identified almost surely as the number of samples grows. \cite{gao2019selecting} considered the case where the covariates only take discrete values and designed an OCBA-based sampling policy that converges to the asymptotic optimal budget allocation rule.

\section{Important research questions on R\&S}
In this section, we outline six R\&S problems that we think are important and yet to be solved. We will explain why we believe these problems are important. Some of these problems have also been considered in the literature. However, we feel that they deserve more research attention.

\vspace{10pt}
\textbf{Problem 1}: Besides the knockout-tournament procedures and the median-elimination procedures introduced in Section 5.2, are there other types of rate-optimal fixed-precision large-scale R\&S procedures?

\textbf{Reason}: Large-scale R\&S is at the center stage of today’s R\&S research because small-scale problems have been studied extensively in the literature, and moreover, they are typically easy to solve with today’s computing resources. The sample-size optimality result of \cite{ZhongHong2019} shows that large-scale problems are fundamentally different from small-scale problems, and many R\&S procedures that are efficient for small-scale problems are not efficient for large-scale problems. Therefore, more procedures need to be proposed under different parallel computing frameworks to solve various large-scale R\&S problems.

\vspace{10pt}

\textbf{Problem 2}: Are there rate-optimal fixed-budget large-scale R\&S procedures?

\textbf{Reason}: Fixed-budget R\&S procedures typically show that their sample-allocation scheme converges to the optimal scheme of \cite{glynn2004large}, as shown in Section 4. However, the optimal scheme
depends heavily on the asymptotic regime, i.e., the number of alternatives stays the same and the total sample size reaches infinity. If the number of alternatives also goes to infinity, as in \cite{ZhongHong2019}, it is no longer optimal. Indeed, the optimal rate of \cite{ZhongHong2019} also applies to fixed-budget R\&S problems. Therefore, fixed-budget large-scale R\&S procedures that are both rate-optimal and efficient in solving practical large-scale problems must be designed.

\vspace{10pt}

\textbf{Problem 3}: How can we design effective dynamic-programming-based procedures to solve fixed-budget R\&S problems?

Reason: As we have shown in Section 4, fixed-budget R\&S problems are in essence finite-time stochastic dynamic programs. However, procedures in the literature are primarily static-allocation approximations or one-step-look-ahead approximations. There appears no serious attempt to directly solve the dynamic programs. However, under Bayesian formulation, the posterior distributions are Gaussian distributions which can be simulated very easily. Therefore, Monte-Carlo-simulation-based approximate dynamic programming (ADP) or reinforcement learning techniques that consider multiple steps seem applicable. Of course, one also has to demonstrate or quantify both theoretically and numerically that going beyond a single step may bring actual benefit.

\vspace{10pt}

\textbf{Problem 4}: How can we design fixed-budget R\&S procedures that are suitable for parallel computing environments?

\textbf{Reason}: As of now, the research attention on parallel R\&S seems primarily focused on fixed-precision procedures. Fixed-budget procedures are often based on dynamic-programming formulation which requires a significant amount of communications among alternatives to determine a sample-allocation policy. Resourceful approaches need be proposed to avoid excessive synchronizations and communications in order to implement fixed-budget procedures in parallel computing environments.

\vspace{10pt}
\textbf{Problem 5}: Do many fixed-precision elimination procedures (e.g., Paulson’s, $\mathcal{KN}$ and $\mathcal{KT}${)} that satisfy PCS guarantee also satisfy PGS guarantee?

\textbf{Reason}: As we reviewed in Section 2, the PCS guarantee requires only a single best, which must be at least $\delta$ larger than all other alternatives. Determining whether a practical problem satisfies this requirement is generally very difficult. Hence, a PGS guarantee that selects an alternative in the indifference zone is certainly more reasonable. However, several fixed-precision elimination procedures (e.g., Paulson’s, $\mathcal{KN}$ and $\mathcal{KT}$) that satisfy the PCS guarantee cannot be proven to satisfy the PGS guarantee. Some of them may be adjusted to satisfy such guarantee at the cost of significantly larger sample sizes \citep{eckman2018guarantees}. To the best of our knowledge, empirical evidence has shown that Paulson’s, $\mathcal{KN}$ and $\mathcal{KT}$ always satisfy the PGS guarantee. Hence, we wonder whether they actually satisfy the PGS guarantee or at least do under some conditions, e.g., when the number of alternatives is large.

\vspace{10pt}
\textbf{Problem 6}: How to better integrate R\&S into optimization-via-simulation (OvS) algorithms?

\textbf{Reason}: Many OvS algorithms require either keeping the current sample best solution or selecting from a group of neighboring solutions. These requirements are naturally R\&S problems. Indeed, a number of R\&S procedures have been proposed to work with OvS algorithms. For instance, \cite{boesel2003using} proposed using R\&S at the end of the OvS process to select the best from all visited solutions, which they called “clean-up”. \cite{pichitlamken2006a} propose to use R\&S in neighborhood selection. Meanwhile, \cite{hong2007a} considered methods to ensure that the current best discovered by the OvS algorithm is indeed the best at any time of the OvS process. In addition, {\cite{he2010simulation} and  \cite{zhang2016optimal} incorporated the OBCA procedures into the cross-entropy and particle swarm OvS algorithms, respectively. However, besides the clean-up procedure, which requires extra work after the OvS process is done and provides no information for OvS algorithms to find a better solution, other ideas tend to slow down the optimization process significantly and output considerably worse solutions. Therefore, it is interesting to figure out how to integrate R\&S into simulation optimization algorithms so that the optimization process may benefit from R\&S. Recently, \cite{10.1145/3170503} showed that the reuse of the simulation observations from the optimization process (also called the search process) by a clean-up procedure at the end of the optimization process may jeopardize the PCS/PGS guarantee.
However, without reusing the simulation observations, R\&S-based clean-up procedures may require a significant amount of simulation observations that may be too expensive to conduct in practice. Therefore, resolving this issue also imposes a theoretical challenge when integrating R\&S into OvS algorithms.}

%\vspace{10pt}
%{\textbf{Problem 7}: How to build the strong connections between  R\&S and MAB problems?}

%{\textbf{Reason}: Even though we have not reviewed R\&S procedures inspired by MAB in details, we think there are a lot of research potentials in building more connections between these two fields. Recent works such as \cite{shin2018tractable} and \cite{glynn2018selecting} bring the ideas from MAB literature to design R\&S procedures in order to achieve
% the optimal decrease rate of PFS derived based on the large deviation results in \cite{glynn2004large}.   When involving the needs of quantifying the input uncertainty,  \cite{WuZhou2019} design R\&S procedures under both the fixed-precision and fixed-budget formulations, and show the convergence rate of the fixed-budget procedure is exponential. How to design R\&S procedures with other rate-optimal measurements in different meaningful asymptotic regimes is worthwhile being investigated.  In addition, it is also interesting to conduct fair comparisons among  different R\&S and MAB algorithms in various practical  problem settings.}

Naturally, the R\&S field has numerous other interesting research problems and emerging research topics. As computing resources are becoming widely available, in general, increasingly complicated R\&S problems need to be solved.

\section*{Acknowledgments}
The authors would like to thank Prof. Shane Henderson for his insightful and detailed comments that have significantly improved this paper. This research was supported in part by the National Natural Science Foundation of China [Grants. 71991473, 71701196, 71722006, and 72031006].

%this research of the second author was supported in part by the Natural Science Foundation of China [Grant 71701196], and this research of the third author was supported in part by the Natural Science Foundation of China [Grants 71722006 and 72031006].
\bibliographystyle{chicago}
\bibliography{RS}
\end{document}